\documentclass{article}

\usepackage[utf8]{inputenc}
\usepackage{graphicx}
\usepackage{subcaption}
\usepackage{amsfonts, amsmath, amssymb, amsthm, authblk}
\usepackage{cleveref}
\usepackage[margin=1in]{geometry}
\usepackage[toc,page]{appendix}
\usepackage{parskip}
\usepackage{url}
\usepackage{tikz}
\usetikzlibrary{arrows.meta,positioning,calc}
\usepackage{booktabs}

\definecolor{colCount}{HTML}{5F5E5A}
\definecolor{colBit}{HTML}{0072b2}
\definecolor{colDest}{HTML}{009e73}

\newtheorem{theorem}{Theorem}[section]

\newtheorem{definition}[theorem]{Definition}

\theoremstyle{remark}
\newtheorem{remark}[theorem]{Remark}

\usepackage{xcolor}

\title{Using Transformers to Generate 3-Manifold Triangulations}
\author[1]{Daniel Bruwel}
\author[2]{Keeley Hoek}
\author[3]{Ashvni Narayanan}
\author[1]{Jonathan Spreer}
\date{}

\affil[1]{{\small School of Mathematics and Statistics F07, The University of Sydney, NSW 2006 Australia, \texttt{daniel.bruwel@sydney.edu.au}, \texttt{jonathan.spreer@sydney.edu.au}}}
\affil[2]{{\small Sydney Mathematical Research Institute, The University of Sydney, NSW 2006 Australia, \texttt{keeley.hoek@sydney.edu.au}}}
\affil[3]{{\small Department of Computer Science, Technion, Haifa, 3200003, Israel, \texttt{ashvni.n@campus.technion.ac.il}}}

\begin{document}

\maketitle

\begin{abstract}
    Using dehydration signatures as an encoding, we train GPT-2-style transformers to generate triangulations of 3-dimensional manifolds. Our main runs each use $\approx 3$ exaFLOPs, spanning a 200 million parameter model trained on 60 million sphere triangulations to an 11 million parameter model trained on 1.1 billion sphere triangulations. The optimal model has 10 blocks, width 640 and generates dehydration signatures of closed manifolds of various topological types at a rate of $6.9\%$, $82\%$ of which are spheres. To investigate how our method scales, we use the technique of IsoFLOP curves on $15$-tetrahedron triangulations of the $3$-sphere, over a compute range of $6$--$3000$ petaFLOPs, and observe neural scaling behaviour. We extract a parameter scaling exponent of $\alpha = 0.49$ ($[0.48, 0.50]$ stat, $\pm 0.02$ syst), and an optimal ratio of $234$ tokens per parameter for training. In addition, we present a number of smaller experiments, including evidence that transformers are able to learn to count degree-one edges near-perfectly, but get progressively worse when counting edges with higher degrees.
\end{abstract}




\section{Introduction}

The study of manifolds is central to research in geometric topology. One way to achieve this is to decompose a given manifold $M$ into simplices, that is, to construct a triangulation $\mathcal T$ of $M$, and to then use the finite combinatorial nature of $\mathcal T$ to compute topological properties of $M$ using algorithms and mathematical software.

For this pipeline of mathematical research to be successful, having access to large numbers of example triangulations with particular properties (determined by the application) is crucial. Ideally, we would have some version of a random generator for triangulations of manifolds that would sample from a controlled distribution over the combinatorial isomorphism types of triangulations with a given number of simplices. The motivation for this requirement stems from algorithms in the field whose performance only depends on the combinatorial isomorphism type of the input, and nothing else.

For manifolds of dimension two, i.e., surfaces, producing such a random generator is well within reach. There exist many random models, see, for instance \cite{Caraceni23Random,Duchon04Random,Pippenger06Random,Poulalhon06Random}, including one by Shankar and the fourth author that is asymptotically uniform on the set of combinatorial isomorphism types of triangulations under consideration \cite{Shankar23RandomSurface}, and one by Altmann and the fourth author based on a Markov chain that converges to a uniform distribution \cite{Altmann26MCMC}.

In dimension $d= 3$, a truly random generator with control over its distribution is much more difficult to obtain. See \cite{Dunfield06Random} for seminal work on random 3-manifolds, including a discussion of random triangulations. Their main investigations, however, focus on a different probabilistic model based on random Heegaard splittings.

If we are interested in sampling triangulations of a given fixed $3$-manifold $M$, we can apply repeated local modifications, e.g., bistellar flips, to a seed triangulation of $M$. Organising this random walk through the space of triangulations of $M$ in a Markov Chain Monte Carlo model achieves asymptotic control over a sampling distribution. See, for instance, \cite{Altmann26MCMC} where this has been done in the case where $M$ is the $3$-sphere. Although convergence rates and mixing times remain unknown, experimental evidence suggests that such methods approximate a random generator of moderately sized triangulations of a given $3$-manifold relatively well.

If we do not want to fix the topology of the triangulation, one straightforward approach is to randomly sample a $4$-regular multi-graph $G = (V,E)$, and for each node $v \in V$, introduce a tetrahedron $\Delta_v$, and for each edge $e \in E$ connecting a vertex $v$ to $w \in V$ (possibly $v=w$), glue a triangle of $\Delta_v$ to a triangle of $\Delta_w$ along a random bijection between their vertices. The issue with this approach is that the probability that this triangulation represents a manifold tends to zero extremely quickly as $\# V \to \infty$, see \Cref{sec:baseline} for some specific numbers.

In this work, we address this issue using machine learning techniques. More specifically, we train a transformer to learn to randomly generate dehydration signatures of $3$-manifold triangulations.

\paragraph{Our contributions} We use GPT-2-style transformers trained on large quantities of dehydration signatures (first presented in \cite{CHW99DehydrationSignature}) of $3$-sphere triangulations. Our hypothesis is that as compute increases, the transformer will improve primarily by increasing the rate of closed manifolds, and to a lesser extent by learning to generate the correct manifold type. This is motivated by the observation that checking if a $3$-dimensional triangulation is a manifold is a linear time procedure, while there is still no polynomial algorithm known for recognising the $3$-sphere \cite{Joswig22SphereRec,Lackenby22Sphere}, arguably the simplest closed $3$-manifold.

We run a series of small experiments on classifying triangulations as manifolds and learning certain properties from their input signatures. We find a large degree of success for simple properties, such as the number of loops in the dual graph, and classification of the manifold property over a small range of topologies. For more complex properties, such as the edge-degree variance, success is more limited.

We run a series of larger experiments for generation. Our models succeed in generating valid triangulations at very high rates, and for triangulations of closed manifolds, at rates many orders of magnitude greater than the baseline. We also find evidence that the optimal parameter count as a function of compute follows a power law, and the optimal loss as a function of compute follows a power law with constant offset.

We observe that our model is significantly biased towards generating spheres, but interestingly that the fraction of generated closed manifolds that are not spheres has a relatively small spread over the range of models we test. In particular, we find it to be only $\approx 1\%$ relative to the spread in the total sphere shortfall. This supports our initial hypothesis that learning is focused on improving the rate of closed manifolds, rather than ensuring those manifolds are the correct type.

Finally, we analyse the statistical properties of our model, and find that there are sphere triangulations with very low likelihoods, preventing our models from being used with techniques such as importance sampling.

\paragraph{Related work} There exist rich and extremely useful data sources for example triangulations, such as censuses of several billion triangulations of manifolds of various topological types in dimensions 3~\cite{BurtonData} and 4~\cite{Burke4Mflds}. By their very nature, these data sources only contain examples with relatively few tetrahedra. To obtain large quantities of triangulations with many tetrahedra ($\geq 12)$, local searches, such as the one from \cite{Altmann26MCMC}, are very useful. However, these only produce a fixed manifold type.

Generating large numbers of example triangulations with distinct topologies can be done using layered handlebodies and random layerings \cite{He25Layerings}, or by operating on framed links \cite{LucyPhD}. Here the combinatorial structure of the examples is very restricted, which is a problem that can be addressed by, again, using local search methods. 

Regarding machine learning on triangulations, in \cite{Costantino_2025}, the authors work with a similar data source, namely, isomorphism signatures of triangulations, to train classifiers to discriminate between topological types, and specific topological properties of the underlying manifold of triangulations. In their work, the machine learning models must pick up on topological features of $3$-manifolds, which are often hard to compute. In contrast, we want our models to learn the manifold property in dimension three, which is a linear time procedure (essentially an Euler characteristic computation).

In a more recent study~\cite{schmidt2026triangulation}, the authors operate on a census of simplicial triangulations of surfaces and $3$-manifolds due to Lutz~\cite{LutzManifoldPage}. They find that the representation of data is very important for learning success, and that performance is highly dependent on the combinatorial structure of the input, along with limited generalisability to triangulations of different sizes.

\paragraph{Acknowledgements}

The authors would like to thank the MATRIX Institute, and in particular the organisers of the MATRIX-Oberwolfach Tandem Workshop {\em Machine Learning and AI for Mathematics}, where this project was started.

The third author would like to thank the Hooper-Shaw Foundation and Sydney Mathematical Research Institute for supporting her as a postdoc when the project began. 
The fourth author is supported by the Australian Research Council under the Discovery Scheme, grant number DP220102588. This research was undertaken with the assistance of resources and services from the National Computational Infrastructure (NCI), supported by the Australian Government, under project {\em io00: Learning triangulations of manifolds}. The authors acknowledge the technical assistance provided by the Sydney Informatics Hub, a Core Research Facility of the University of Sydney. 

\paragraph{AI disclosure statement:} The authors acknowledge the use of generative AI---specifically, web-based ChatGPT 6.0 and Claude Opus 5.5---for proofreading, for writing some of the TikZ code, and for assistance with minor coding tasks. 

\section{Background}

\subsection{Manifolds}

A {\em manifold} of dimension $d$ is a space that locally looks like Euclidean space $\mathbb{R}^d$, but globally may have very different properties. The prime example of a manifold is the $2$-dimensional sphere, that is, the surface of the Earth. For an example in arbitrary dimension $d$, consider the $d$-fold Cartesian product of the circle, which, for $d=2$, forms the surface of a doughnut. For an interesting family of $3$-dimensional manifolds, consider two solid tori, glued along their boundaries such that a closed curve bounding a disc (called the {\em meridian}) on one side is glued to an arbitrarily complicated simple closed curve on the other side. This forms an infinite family of manifolds called {\em lens spaces}, parameterised by two coprime integers $p > 0$ and $q \geq 0$ and denoted by $L(p,q)$. 
The family contains the $3$-dimensional sphere, that is, the space of points in $\mathbb{R}^4$ at unit distance from the origin.

Manifolds play an important role in many fields of mathematics, and other sciences such as physics or chemistry, or disciplines in engineering such as robotics.
We consider manifolds from the viewpoint of geometric topology. That is, we see two manifolds as equivalent, or, {\em homeomorphic}, if one can be transformed into the other by a continuous deformation. Specifically, we focus on manifolds in dimension $d=3$. This is an interesting special case for many reasons. First, in this case, every topological 3-manifold has a unique piecewise linear (PL) structure. This means that the topological (continuous deformations) and PL (piecewise linear deformations) classifications agree, making this dimension very suitable for studying manifolds using triangulations (see below). Second, dimension $d=3$ is the last dimension where we have an algorithmic solution to the {\em homeomorphism problem}, that is, the task of deciding whether two given manifolds are equivalent \cite{AschenbrennerFriedlWilton2015, Markov}.

See \cite{hatcher, jaco} for more background on 3-manifolds, and \cite{armstrong} for a general overview of basic topology.

\subsection{Triangulations}

A sufficiently well-behaved manifold, specifically one that is {\em piecewise linear (PL)}, admits a decomposition into finitely many simplices. To this end and in dimension $d=3$, consider $N_{\text{tet}}$ disjoint tetrahedra $\Delta_1, \ldots , \Delta_{N_{\text{tet}}}$ together with a set of $2 N_{\text{tet}}$ affine maps $\Phi_{1}, \ldots , \Phi_{2 N_{\text{tet}}}$ identifying their $4 N_{\text{tet}}$ boundary triangles in pairs. Then the quotient $\mathcal{T} = \{ \Delta_1, \ldots , \Delta_{N_{\text{tet}}} \} / \{ \Phi_{1}, \ldots , \Phi_{2 N_{\text{tet}}} \}$ forms an underlying $3$-dimensional space, denoted by $|\mathcal{T}|$ that can be thought of as the underlying space obtained from gluing $N_{\text{tet}}$ tetrahedra together along their triangular faces (vertex to vertex).

We refer to such a $\mathcal{T}$ as a {\em triangulation}. As a result of the face gluings, many vertices, edges, and triangles of the $\Delta_i$, $ 1\leq i \leq N_{\text{tet}}$, become identified in $\mathcal{T}$ and we will refer to them as vertices, edges and triangles of the triangulation $\mathcal{T}$. We refer to the number of tetrahedra $N_{\text{tet}}$ in a triangulation as its {\em size}.

We say that $\mathcal{T}$ is {\em valid}, if no edge of any of the tetrahedra becomes identified with itself in reverse as a consequence of the face gluings. Moreover, we say that $\mathcal{T}$ is a {\em triangulated manifold} or {\em closed} if, in addition, its {\em vertex links}, that is the boundaries of small neighbourhoods of the vertices of $\mathcal{T}$, are all $2$-dimensional spheres. Note that this is equivalent to checking that all points in $|\mathcal{T}|$ locally look like Euclidean space. Finally, given a manifold $\mathcal{M}$, we say that $\mathcal{T}$ {\em triangulates} $\mathcal{M}$ if $|\mathcal{T}|$ is homeomorphic to $\mathcal{M}$.

Given a triangulation $\mathcal{T}$, its {\em Euler characteristic} is given by the alternating sum over the number of faces in every dimension. For a $3$-dimensional triangulation this is the number of vertices minus the number of edges plus the number of triangles minus the number of tetrahedra. The Euler characteristic is a {\em topological invariant}. This means that, given a manifold $\mathcal{M}$, every triangulation of $\mathcal{M}$ must have the same Euler characteristic, but two triangulations with equal Euler characteristics may triangulate distinct manifolds. In dimension three (in fact, in all odd dimensions), all closed manifolds must have Euler characteristic $0$. In contrast, in dimension $2$, the topological type of a connected, closed, orientable surface is determined by its Euler characteristic \cite{GallierXu2013}. 

If the link of a vertex in a valid triangulation is a surface other than the sphere, we call this vertex {\em ideal}, and think of a small open neighbourhood of it as being removed from the triangulation. This implicitly creates a boundary component of the underlying space of $\mathcal{T}$. In this case we call $\mathcal{T}$ an {\em ideal triangulation}\footnote{We note here that we do not require all vertices of an ideal triangulation to be ideal. This is somewhat in contrast to the literature on $3$-manifold triangulations.}. Given an ideal triangulation, we will use the Euler characteristic of its vertex links as a heuristic diagnostic for how much the triangulation fails to be closed. Here, the ``distance'' from being a closed manifold is captured by how far away this quantity is from the Euler characteristic of the $2$-sphere.

\subsection{Dehydration signatures}
\label{ssec:dehydration}

A dehydration signature, as first described in \cite[Section 4]{CHW99DehydrationSignature}, is an encoding of a triangulation with $N_{\text{tet}} \leq 25$ tetrahedra into a string of $\approx 5 N_{\text{tet}}/2$ lowercase English alphabet characters. It consists of four sections: a ``size'' character, a run of ``new tetrahedra'' characters, $N_{\text{tet}}+1$ ``where'' characters, and $N_{\text{tet}}+1$ ``how'' characters. 

The ``size'' character encodes the number of tetrahedra that the triangulation has: $a$ means the triangulation is empty, $b$ means it has a single tetrahedron and so on. The ``new tetrahedra'' characters encode ordinal information about each tetrahedron. Specifically, while looping through all $2 N_{\text{tet}}$ face pairings, they encode when new tetrahedra are first introduced. For this, the signature reserves $2 \lfloor (N_{\text{tet}}+3)/ 4 \rfloor$ characters that are read in pairs $XY$. Each $X$ and $Y$ takes one of $16$ letters letting the pair encode a bit string of length $8$, encoding which of $8$ gluings add a new tetrahedron. Note, these tetrahedra are glued together along identity permutations, which gives sufficient information to glue all $N_{\text{tet}}$ tetrahedra along a spanning tree of the dual graph of the triangulation, and to fix the labelling. This takes care of $N_{\text{tet}}-1$ of the overall $2 N_{\text{tet}}$ face pairings.

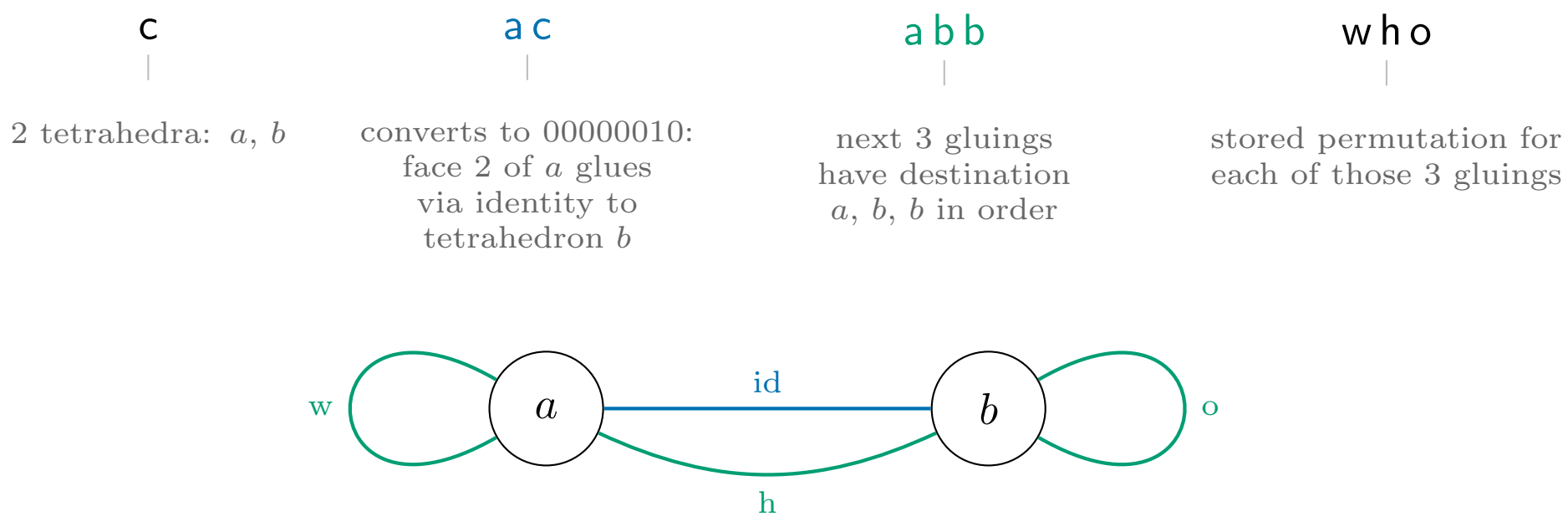
\begin{figure}[htpb]
    \centering
    \begin{tikzpicture}[
  font=\sffamily,
  desc/.style={font=\scriptsize, text=black!60, align=center, text width=2.8cm},
  tick/.style={gray!50, thin}
]
\node (c)   at (0,0)    {c};
\node (ac)  at (3.0,0)  {\textcolor{colBit}{a\,c}};
\node (abb) at (6.3,0)  {\textcolor{colDest}{a\,b\,b}};
\node (who) at (9.8,0)  {w\,h\,o};

\draw[tick] (c.south) -- ++(0,-2mm);
\node[desc, below=4mm of c] {2 tetrahedra: $a$, $b$};

\draw[tick] (ac.south) -- ++(0,-2mm);
\node[desc, below=4mm of ac] {converts to 00000010: face 2 of $a$ glues via identity to tetrahedron $b$};

\draw[tick] (abb.south) -- ++(0,-2mm);
\node[desc, below=4mm of abb] {next 3 gluings have destination $a$, $b$, $b$ in order};

\draw[tick] (who.south) -- ++(0,-2mm);
\node[desc, below=4mm of who] {stored permutation for each of those 3 gluings};

\begin{scope}[xshift=3.15cm, yshift=-3cm]
\node[draw, circle, minimum size=9mm] (a) at (0,0)   {$a$};
\node[draw, circle, minimum size=9mm] (b) at (3.5,0) {$b$};

\draw[colBit, thick] (a) -- (b)
  node[midway, above, font=\scriptsize] {id};

\draw[colDest, thick] (a) to[out=150,in=210,looseness=10]
  node[left, font=\scriptsize] {w} (a);

\draw[colDest, thick] (a) to[bend right=25]
  node[below, font=\scriptsize] {h} (b);

\draw[colDest, thick] (b) to[out=30,in=-30,looseness=10]
  node[right, font=\scriptsize] {o} (b);
\end{scope}
\end{tikzpicture}
    \caption{Example dehydration signature for a single vertex triangulation of $S^3$ --- \texttt{cacabbwho}}
    \label{fig:dehydration-schema}
\end{figure}

To identify the remaining $N_{\text{tet}}+1$ gluings, the next $N_{\text{tet}}+1$ ``where'' characters encode to which tetrahedron the remaining open faces are glued. Here, the order of the open faces is determined by the previous step. For each target tetrahedron, one letter is used to encode an integer between $0$ and $25$. Finally, the last $N_{\text{tet}}+1$ ``how'' characters determine the permutation of the remaining face pairings. Here, note that there are $4! = 24$ possible gluings, since not only the vertices of a triangle must be glued, but also the face of the target tetrahedron must be identified. For this, one character with possible values ``$a$'' encoding $0123$ to ``$x$'' encoding $3210$ is used. These permutations are in ``reverse bijection notation''. I.e. $0123$ maps $3\to 0$, $2\to 1$, $1\to 2$, $0\to 3$ and $3210$ is the identity. \Cref{fig:dehydration-schema} demonstrates this gluing for a single vertex triangulation of $S^3$. All indices start at $0$.

\begin{remark}
    We use the term \textit{dehydration signature} in a more general sense than other literature. In particular, we mean any signature that can be decoded by the above technique and produce valid gluing information, even if it is not strictly the dehydration of a triangulation. See \Cref{sec:baseline} for more details.
\end{remark}
\begin{remark}
Note that this encoding only slightly masks the combinatorics of the underlying triangulation. This is in contrast to more modern and compact encodings such as the isomorphism signatures used by {\em Regina} \cite{regina}. It is reasonable to assume that this has a positive effect on learning by transformers, as suggested by preliminary tests.
\end{remark}
\begin{remark}
For similar experiments involving triangulations with more than $25$ tetrahedra a slight tweak to the dehydration signature format is necessary. This, as well as investigating other machine-learning-friendly triangulation encodings, is subject of future work.    
\end{remark}

\section{Baseline Generation Efficiency}\label{sec:baseline}
Before commencing our examination of deep learning techniques, we establish a series of baselines. Our investigations are done entirely on generating \textit{dehydration signatures}, see \Cref{ssec:dehydration} for details.
We note that the generation technique we will consider is autoregressive, that is to say it generates the dehydration signatures character by character. This motivates the following definitions.

\begin{definition}~
    \begin{itemize}
        \item A string is \textbf{acceptable} if it can be parsed as a dehydration signature
        \item A dehydration signature is \textbf{valid} if the gluing rules encode a valid triangulation. That is, the gluing rules do not identify an edge with itself in reverse.
        \item A dehydration signature is \textbf{closed} if the triangulation it represents is a valid closed manifold. That is, the link of each vertex is PL-homeomorphic to $S^2$.
        \item A dehydration signature is a \textbf{sphere} if the triangulation it represents is PL-homeomorphic to an $S^3$.
    \end{itemize}
\end{definition}

By construction, we have
\[\text{sphere}\subset\text{closed}\subset\text{valid}\subset\text{acceptable}\]

\paragraph{Acceptable Dehydration Signatures}
As a naive baseline, we want to establish the probability that a random sequence of lower case English letters, of the correct length, forms an acceptable dehydration signature. Exact enumeration of dehydration signatures is challenging, primarily due to the fact that ``forward references'' are possible: that is, a face can be glued to a tetrahedron that has not yet been defined to exist in the ``new tetrahedra'' characters of the string. Instead of exact enumeration, we use probabilistic estimates on the tree of acceptable dehydration signatures, as introduced in \cite{knuth1975estimating}. Details are provided in \Cref{appendix:estimating-dehydration-count}. \Cref{tab:baseline} includes our estimation of the probabilities: column 2 contains $P(\text{Acceptable} \mid \text{Random})$, the probabilities when sampling each character unrestricted from all $26$ lower case English letters; column 3 contains $P(\text{Acceptable} \mid \text{Subset})$, the probabilities when sampling each character from a restricted ``subset'' of available characters at each position (e.g. the letter \texttt{k} is always required for the first position of a 10-tetrahedron signature).

We take a moment to define an important subclass of acceptable dehydration signatures called \textbf{encodable} dehydration signatures. These are defined as acceptable dehydration signatures that result from the encoding of some triangulation via the standard dehydration algorithm \cite{CHW99DehydrationSignature}. These are important because all training data will only consist of encodable dehydration signatures. For a fixed triangulation, the number of unique encodable dehydration signatures, $s(\mathcal{T})$, is given by

\begin{equation}\label{eq:sig-count}
    s(\mathcal{T})=\frac{24N_{\text{tet}}(\mathcal{T})}{\#\operatorname{Aut}(\mathcal{T})}
\end{equation}

where $N_{\text{tet}}(\mathcal{T})$ denotes the number of tetrahedra in $\mathcal{T}$, and $\operatorname{Aut}(\mathcal{T})$ denotes its group of symmetries.

\paragraph{Valid Triangulations}
We now consider the conditional probability of a triangulation being valid, given that it has an acceptable dehydration signature. We estimate this by randomly sampling from the space of triangulations directly. We first note an important detail: a triangulation has a dehydration signature if and only if it is connected and has no boundary faces. To generate such triangulations we randomly pair faces together to form a perfect matching, and assign a random permutation to each pairing. This is guaranteed to form a triangulation with no boundary faces. Our sampled triangulation may not be connected so we discard all disconnected triangulations. This is efficient as asymptotically almost surely a randomly generated 4-valent graph is connected, and for $N_{\text{tet}}=10$ the probability of being connected is $>95\%$. We generate 100,000 such triangulations for each $N_{\text{tet}} = 10$ to $20$ and find that $P(\text{Valid}\;|\;\text{Acceptable})$, the probability of being valid given an acceptable signature, is on the order of 14-18\%. This is included in the fourth column of \Cref{tab:baseline}.

\paragraph{Triangulations of Closed Manifolds}
For $N_{\text{tet}}=5, 6$ we use the above random gluing method to generate connected triangulations with no boundary faces. For each $N_{\text{tet}}$ we take 5 million random gluings (of which in both cases $\approx$ 1 million are valid) and find that $P(\text{Closed} \;|\;\text{Acceptable})$, the probability of being closed given the signature is acceptable, for $N_{\text{tet}}=5$ is 0.053\% and for $N_{\text{tet}}=6$ is 0.005\%. For $N_{\text{tet}}=10$ we test 1 billion samples and find not a single triangulation of a closed manifold. Hence, with 95\% confidence we can say that $P(\text{Closed} \;|\;\text{Acceptable})$ is $<3\times 10^{-9}$. It is known that the probability of being closed tends to $0$ as the number of tetrahedra tends to infinity \cite[Proposition 2.8]{Dunfield06Random}. Hence, direct sampling with the above method can be assumed to be inefficient for $N_{\text{tet}}\geq 10$.

\paragraph{Triangulations of Spheres}
Considering the above manifold generation is inefficient for generating closed manifolds, we cannot use it to estimate the percentage of closed manifolds that are spheres. We also cannot practically use the other random sampling techniques such as layered handlebodies or random layerings as the combinatorial type is not uniform. Hence, we look to the exhaustive census of closed triangulations of manifolds for $1 \leq N_{\text{tet}} \leq 9$ (with arbitrary numbers of vertices) \cite{Burton11Census}.
From this, approximately 31\% of closed triangulations are spheres, the breakdown for $1 \leq N_{\text{tet}}\leq 9$ is shown in \Cref{tab:sphere-rate}. For larger $N_{\text{tet}}$ the census is incomplete, though it is reasonable to assume that, as $N_{\text{tet}}$ increases, the percentage of closed manifolds that are spheres declines.

\begin{table}[htpb]
    \centering
    \begin{tabular}{c|c|c|c}
        $N_{\text{tet}}$ & $P(\text{Acceptable} \mid \text{Random})$ & $P(\text{Acceptable} \mid \text{Subset})$ & $P(\text{Valid} \mid \text{Acceptable})$ ($\%$)  \\
        \hline
        $10$ & $1.86\times 10^{-15}$ & $7.88\times 10^{-8}$ & $17.70\%$ \\
        $11$ & $5.23\times 10^{-16}$ & $1.99\times 10^{-8}$ & $17.40\%$ \\
        $12$ & $1.71\times 10^{-16}$ & $5.37\times 10^{-9}$ & $17.00\%$ \\
        $13$ & $2.65\times 10^{-17}$ & $1.68\times 10^{-9}$ & $16.60\%$ \\
        $14$ & $9.53\times 10^{-18}$ & $4.31\times 10^{-10}$ & $16.30\%$ \\
        $15$ & $4.00\times 10^{-18}$ & $1.21\times 10^{-10}$ & $15.90\%$ \\
        $16$ & $1.66\times 10^{-18}$ & $3.14\times 10^{-11}$ & $15.70\%$ \\
        $17$ & $2.85\times 10^{-19}$ & $8.42\times 10^{-12}$ & $15.50\%$ \\
        $18$ & $2.19\times 10^{-19}$ & $3.62\times 10^{-12}$ & $15.30\%$ \\
        $19$ & $8.47\times 10^{-20}$ & $7.43\times 10^{-13}$ & $15.00\%$ \\
        $20$ & $6.05\times 10^{-20}$ & $2.68\times 10^{-13}$ & $14.80\%$ \\
        
    \end{tabular}
    \caption{Baseline generation efficiency of valid and acceptable dehydration signatures.}
    \label{tab:baseline}
\end{table}

\begin{table}[htpb]
    \centering
    \begin{tabular}{c|c|c|c}
        $N_{\text{tet}}$ & Closed & Sphere & Sphere Percent \\
        \hline
        1 & $4$ & $2$ & $50\%$ \\
        2 & $17$ & $6$ & $35\%$ \\
        3 & $81$ & $32$ & $40\%$ \\
        4 & $577$ & $198$ & $34\%$ \\
        5 & $5{,}184$ & $1{,}903$ & $37\%$ \\
        6 & $57{,}753$ & $19{,}935$ & $35\%$ \\
        7 & $722{,}765$ & $247{,}644$ & $34\%$ \\
        8 & $9{,}787{,}509$ & $3{,}185{,}275$ & $33\%$ \\
        9 & $139{,}103{,}032$ & $43{,}461{,}431$ & $31\%$ \\
        \hline
        \textbf{Total} & $\mathbf{149{,}676{,}922}$ & $\mathbf{46{,}916{,}426}$ & $\mathbf{31\%}$
    \end{tabular}
    \caption{Sphere triangulation rates for small triangulations \cite{Burton11Census}}
    \label{tab:sphere-rate}
\end{table}

\section{Data and Model Description}
\label{sec:data-model-description}
\subsection{Training Data Generation}
\paragraph{Triangulations}
While our primary goal is to study the generation of triangulations, we will first establish a collection of classification and regression tasks. These serve as critical diagnostic tools to understand if transformer architectures can inherently learn to compute fundamental topological and combinatorial properties from dehydration signatures. To support this and the final generation tasks, we generate four sets of triangulation data. We generate most of these data using MCMC as discussed in \cite{Altmann26MCMC}. We use a step size of $d_{\text{step}}=1$ and choose $\gamma$, the growth control factor, according to the number of tetrahedra that we are trying to generate, extracted via the method outlined in Section~4.1 of \cite{Altmann26MCMC}.

First, for each $10 \leq N_{\text{tet}}\leq 20$, we generate 4 million unique triangulations of $S^3$ with one vertex. Second, we generate 1.2 billion unique samples of $S^3$ with $N_{\text{tet}}=15$, with one vertex. Third, we generate $1{,}000$ unique samples for the following manifolds: $L(1,0)=S^3$, $S^2\times S^1$, $L(2,1)=\mathbb{R}P^3$, $L(7, 1)$, and $L(7, 2)$, all with $N_{\text{tet}}=10$. Finally, we generate a single ``manifold learning'' dataset with $N_{\text{tet}}=10$. This dataset consists of
\begin{itemize}
    \item 1 million unique samples of $S^3$ (a subset from above).
    \item 1 million unique triangulations of manifolds generated by randomly selecting a triangulation from the census of ``internal, finite'' triangulations of 6 tetrahedra, and performing 4 random (2-3) moves.
    \item 1 million unique valid, connected triangulations with no boundary faces generated by random gluings and rejection sampling as outlined in \Cref{sec:baseline}.
    \item 1 million unique random ``near miss'' triangulations,  generated by taking a pair of triangulations from the census of ``internal, finite'' triangulations of 6 tetrahedra, removing one tetrahedron from each triangulation (chosen at random among those without self-identified triangles), and then randomly gluing the 4 faces of the first punctured manifold with the 4 faces of the second punctured manifold. This is also run with rejection sampling to ensure that the resulting triangulation is connected, valid, and not a manifold.
\end{itemize}
This gives 4 million unique triangulations all without boundary faces, valid, and connected, half of which are manifolds, and half of which are not.

All the above triangulations are saved with their dehydration signatures, and no canonicalisation is done. Additionally, uniqueness is qualified on the signature level, not the triangulation level; i.e., two combinatorially identical triangulations with different labellings are considered to be distinct.

\paragraph{Scalar Objective Functions}
For each of the 4 million samples of $S^3$ with one vertex and $N_{\text{tet}}=10$, we compute the following objective functions.
\begin{enumerate}
    \item The number of self loops in the dual graph.
    \item The number of edges with degree one, two, three, four, and five.
    \item The variance of the edge degree.
\end{enumerate}

\subsection{Model}

We use transformer models, first described in \cite{GoogleAttention}. More precisely, our model is based on OpenAI's GPT-2 transformer introduced in \cite{radford2019language}.

\label{subsec:model-details}
\paragraph{Signature Tokenisation}
Transformer models are built to expect a list of ``tokens'' as their input. These tokens are discrete elements of some ``vocabulary''. Hence we need to tokenise our dehydration signatures. See \Cref{ssec:dehydration} for a detailed description of its format. Given a dehydration signature, we tokenise each character independently, accounting for the letter (\texttt{a, b, c, ...}) and its role (\texttt{size, new tetrahedra, where, how}). The vocabulary consists of tokens such as \texttt{size <a>, size <b>, where <a>} etc. with an overall vocabulary size of $1$ start of sentence character \texttt{[SOS]}, $1$ end of sentence character \texttt{[EOS]}, $1$ ``size'' character, $16$ ``new tetrahedra'' characters, $N_{\text{tet}}$ ``where'' characters, and $24$ ``how'' characters, totaling $43+N_{\text{tet}}$. 

\paragraph{Core Model Design}
For the core transformer, we use a standard GPT-2-style architecture from \cite{radford2019language}. Here, we briefly review its structure for convenience. The transformer consists of learned positional and token embeddings into $\mathbb{R}^{d_{\text{model}}}$ that are added together. This is followed by a dropout layer, a series of transformer blocks, a layer norm and a final dense layer projecting to $\mathbb{R}^{d_{\text{out}}}$. Each of the transformer blocks consists of an initial layer norm for stability, a residually connected attention layer, a dropout layer, a second layer norm, a residually connected feed forward layer, and a final dropout. The feed forward layer is a simple dense projection into $\mathbb{R}^{4\cdot d_{\text{model}}}$ space, a Gaussian Error Linear Unit (GELU) activation \cite{hendrycks2023gaussianerrorlinearunits}, and a projection back into $\mathbb{R}^{d_{\text{model}}}$. The hyperparameters we need to choose to determine the shape of the model are the embedding dimension $d_{\text{model}}$, the number of transformer blocks $n_{\text{block}}$, the number of attention heads per block $n_{\text{head}}$, and output dimension $d_{\text{out}}$. Additionally, each transformer block can include a causal attention mask.

The total number of model parameters can be approximated as $12 \cdot d_{\text{model}}^2\cdot n_{\text{block}}$. Note that in any case where the parameter count is explicitly stated, the exact parameter count is used rather than this approximation.

\begin{figure}[htpb]
    \centering
    \begin{tikzpicture}[
    grouped/.style={draw, rectangle, minimum width=4.5cm, minimum height=0.7cm, align=center, fill=blue!5, rounded corners},
    pooling/.style={draw, rectangle, minimum width=4.5cm, minimum height=0.7cm, align=center, fill=red!5, rounded corners},
    optional/.style={draw, rectangle, minimum width=4.5cm, minimum height=0.7cm, align=center, fill=green!5, rounded corners, dashed},
    arrow/.style={-Latex, thick},
    scale=1, transform shape
]

\node (input) {$[v_0, \dots, v_n]$};

\node[grouped, below=0.6cm of input] (ln1) {Layer Normalisation};

\node[grouped, below=0.6cm of ln1] (weights) {Dense (1D) $\rightarrow$ Softmax \\ (Compute weights $w_i$)};

\node[pooling, below=0.6cm of weights] (pooling) {Weighted Sum \\ $v = \sum_i w_i \hat{v}_i$};

\node[grouped, below=0.6cm of pooling] (mlp) {Layer Norm $\rightarrow$ Dense\\ $\rightarrow$ ReLU $\rightarrow$ Scalar Projection};

\node[optional, below=0.6cm of mlp] (output) {Sigmoid};
\node[below=0.6cm of output] (loss) {Final Output};

\draw[arrow] (input) -- (ln1);
\draw[arrow] (ln1) -- (weights) node[midway, right] {$[\hat{v}_0, \dots, \hat{v}_n]$};
\draw[arrow] (weights) -- (pooling);
\draw[arrow] (pooling) -- (mlp) node[midway, right] {$v$};
\draw[arrow] (mlp) -- (output);
\draw[arrow] (output) -- (loss);

\draw[arrow, dashed] (ln1.west) -- ++(-1.5,0) |- (pooling.west) node[near start, left] {$[\hat{v}_i]$};

\end{tikzpicture}
    \caption{Architecture for the classification/regression head. The final sigmoid term is only included during classification.}
    \label{fig:regression-head}
\end{figure}
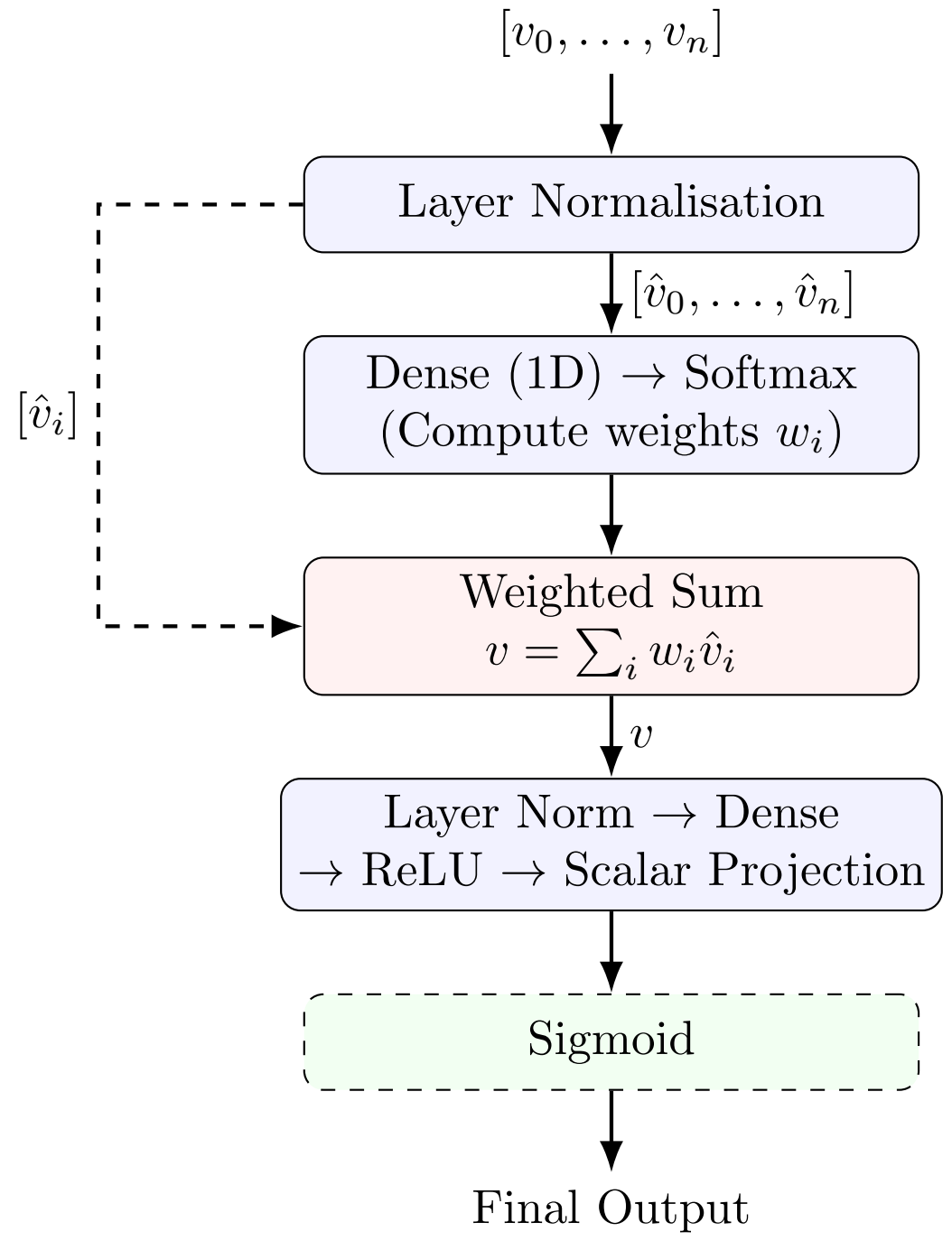

\paragraph{Classification and Scalar Regression Model}
For the scalar regression objectives, we use the main transformer block from above specifying $d_{\text{out}}$. This produces a list of vectors $v_0,\ldots,v_n \in \mathbb{R}^{d_{\text{out}}}$ corresponding to each input token. To convert these to a scalar value, we use a custom scalar head illustrated in \Cref{fig:regression-head}. The scalar head takes each of the output vectors, applies a layer normalisation to them to produce $\hat{v}_0,\ldots,\hat{v}_n \in \mathbb{R}^{d_{\text{out}}}$. These normalised vectors are then projected down to one dimensional logits $l_0,\ldots,l_n \in \mathbb{R}$ which are then passed through a softmax to obtain weights $w_i = \frac{e^{l_i}}{\sum_j e^{l_j}} \in \mathbb{R}$. These are pooled together to form $v=\sum_i w_i\hat{v}_i$.

This is then passed through another layer normalisation, and then a dense layer, a Rectified Linear Unit (ReLU) activation function, and finally projected down to a single scalar value. For the classification models, this scalar value is then passed through a final sigmoid function and the loss function used is binary cross entropy. The scalar regression models are all trained with mean squared error as the loss, and the data is normalised to a mean of 0 and standard deviation of 1 prior to training. 

Additionally, there is some ambiguity whether causal masking should be used or not. On the one hand, the masking can introduce inductive bias and lead to faster convergence, on the other hand it may limit the ability for tokens to attend to each other. Initial testing found that causal masking performed better and is included in all the models.

\paragraph{Autoregression Model}
For autoregression, we take the standard transformer with $d_{\text{out}}$ equal to the size of the vocabulary. We take our encoded signature and prepend a start of sentence token \texttt{[SOS]} to it. The model is trained with the objective of predicting the next token in the sequence. That is, we want to take a sequence \texttt{[SOS, $T_1$, $T_2$, $T_3$, ..., $T_n$]} and, for each token $T_i$, predict $T_{i+1}$, that is, predict \texttt{[$T_1$, $T_2$, $T_3$, ..., $T_n$, EOS]}. The output represents the logits for each token, and standard cross entropy is used. Causal masking is included to ensure the model does not use future tokens to predict the current token.

The dropout rate is fixed at $0.1$ for training across all models. This includes single epoch runs where dropout is not needed, and may have degraded performance slightly.

Generation is done by initialising a sequence to \texttt{[SOS]}, predicting the next token, appending to extend the sequence, and repeating. At each point, the next token is decided by taking the output logits over the vocabulary, applying a softmax to form a list of probabilities, and then randomly sampling according to these probabilities.

Training is done in \texttt{bfloat16}. Small models are trained on NVIDIA A5000 GPUs, and large models are trained on NVIDIA H200 GPUs. 

\paragraph{Training}
All models use the stochastic optimisation algorithm AdamW to reduce the training error. We operate on a batch size of 512. AdamW is chosen for its robustness to hyperparameter choice, and its fast wall-clock convergence. For AdamW there are 4 hyperparameters. We keep the weight decay, $\beta_1$ and $\beta_2$ at the standard values of $0.01$, $0.9$ and $0.999$ respectively for all models. For the learning rate we use a standard learning rate schedule consisting of a linear warm-up over the first $\approx 5\%$ of the iterations to a predetermined peak learning rate, followed immediately by a cosine annealing schedule down to $10\%$ of the peak learning rate \cite{he2018bagtricksimageclassification}. This learning rate schedule is very stable, but because it depends on the total number of iterations, a checkpoint taken partway through a run is not equivalent to a fully trained model with a smaller budget; we therefore train every (model, budget) pair from scratch.

The peak learning rate needs to be determined independently for each model and each objective function to prevent conflation between architectural and objective function performance with suboptimal learning rate choices. That is, we avoid saying that one objective function is easier to learn than another, or one architecture is better than another, if the performance difference is simply due to a poor learning rate choice. 

To choose the peak learning rate we rely on two techniques. First, for all the small tests, we run a coarse grid search of 3-5 learning rates, approximately logarithmically spaced, and pick the optimum. In each case, the models are trained to completion and the optimum is chosen based on a combination of the final test loss and the training loss profile. This strategy is chosen due to its simplicity and to avoid spending large amounts of compute on fine grained hyperparameter tuning where it is not necessary. Second, for large scale tests, we use the full $\mu P$ parametrisation of \cite{yang2022tensorprogramsvtuning}. Concretely, $\mu P$ is a parametrisation under which a \textit{base} learning rate is fixed, and different layers get different learning rate multipliers and initialisation variances depending on their role. In this parametrisation, the optimal \textit{base} learning rate is approximately independent of width ($d_{\text{model}}$), so it can be tuned on a narrow model and transferred.

\section{Numeric Experiments}
\label{sec:experiments}
While the primary goal of this study is to assess the efficiency of transformers for autoregression of dehydration signatures, we include some experiments on classification and regression tasks to assess and diagnose the transformer architecture's performance on learning the combinatorial properties of dehydration signatures.

This is important as it helps diagnose what sorts of combinatorial properties transformers can recognise. Often, simple properties give important insights to the nature of a triangulation. For instance, low-degree edges  often permit simplification, and under standard hypotheses they are excluded from minimal triangulations \cite{Jaco09LensSpaces}. Another example is given by a conjecture due to Thurston and proven by Elder, McCammond and Meier \cite{ELDER20031241}, stating that a triangulation of a manifold with all edges having degree five or six and no two degree five edges in a single triangle must have a word-hyperbolic fundamental group.

Hence, we start this section by looking at the learnability of combinatorial properties of triangulations in \Cref{subsec:manifold-property,subsec:scalar-regression} before moving on to our main experiment in \Cref{subsec:auto-regression}.

\subsection{Detection of the Manifold Property}
\label{subsec:manifold-property}

In this first experiment we want to investigate how reliably transformers
can identify closed triangulations in a larger list of more general 
triangulations.

We run a single model architecture consisting of $d_{\text{model}}=512$, $n_{\text{head}}=8$, $n_{\text{block}}=8$, and $d_{\text{out}} = 64$. This model has $\approx 25$ million parameters. We therefore trained the model for 16 epochs on the 4 million samples in the ``manifold learning'' dataset described in \Cref{sec:data-model-description}, leading to approximately 2.5 samples per parameter, setting aside $1{,}000$ of each of the sphere, random manifold, random gluing, and near miss datasets for testing. We also include the $1{,}000$ samples of $L(7, 1), L(7, 2), \mathbb{R}P^3$ and $S^2\times S^1$ in the testing group for a total of $8,000$ test samples. We use a coarse peak learning rate grid of $3\times 10^{-5}, 1\times10^{-4}, 3\times 10^{-4},$ and $1\times10^{-3}$, determining the optimal value as $1\times10^{-4}$. This is a small model and a small dataset with a high number of epochs. Consequently, the results are limited and would likely be surpassed by larger models with more training data.

\begin{figure}[htbp]
    \centering
    \begin{subfigure}[b]{0.48\linewidth}
        \centering
        \includegraphics[width=\linewidth]{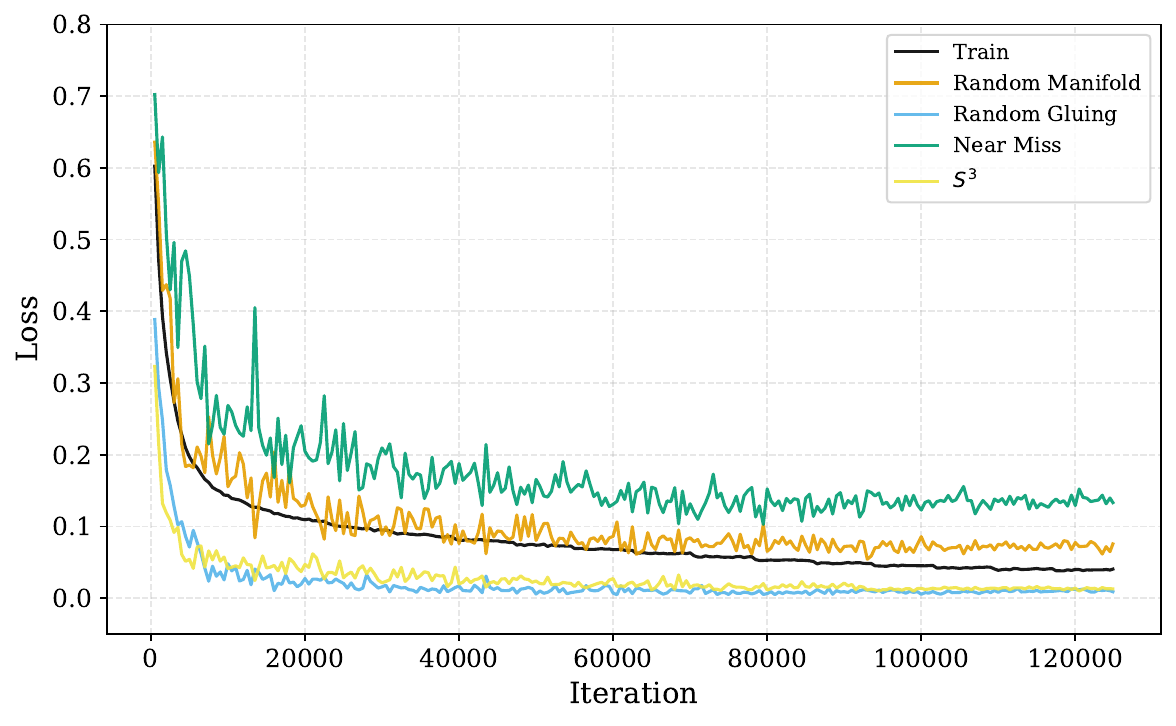}
        \caption{}
    \end{subfigure}
    \hfill
    \begin{subfigure}[b]{0.48\linewidth}
        \centering
        \includegraphics[width=\linewidth]{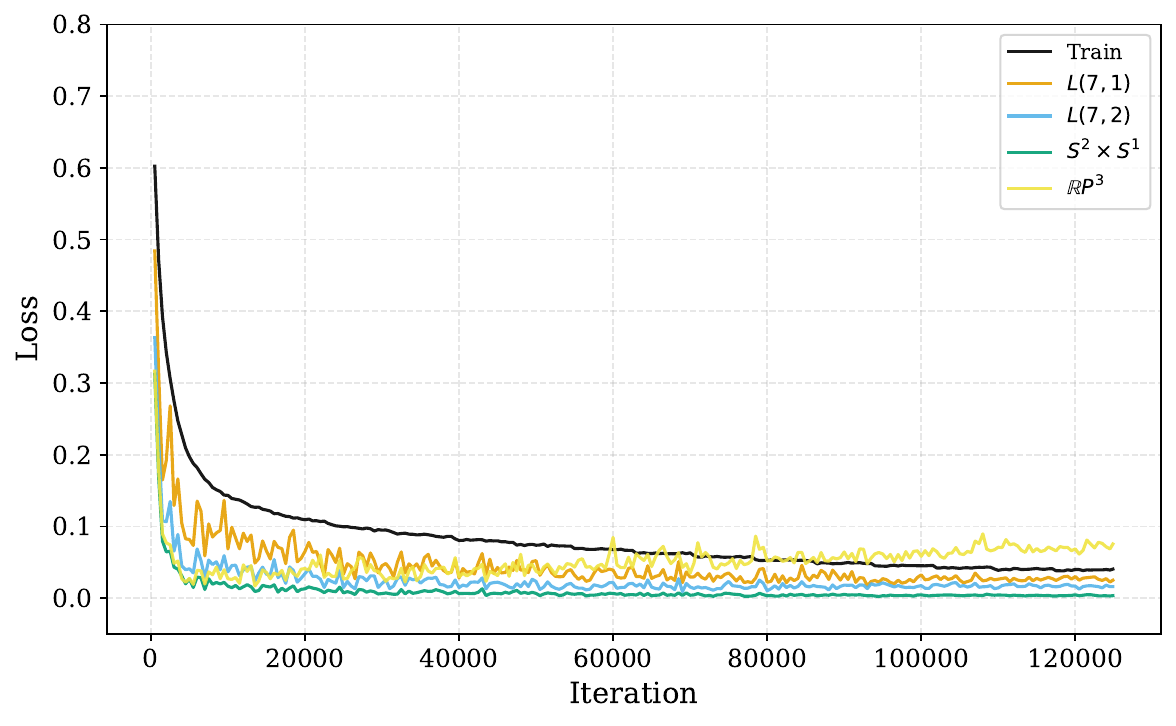}
        \caption{}
    \end{subfigure}
    \caption{Loss curves for manifold classification task.}
    \label{fig:classification-loss}
\end{figure}

The loss profile between the different manifold classes and the training loss is shown in \Cref{fig:classification-loss}. The loss for all manifold types follows a consistent profile over the duration of the training. However, there is a substantial difference in the scale of the loss. We extract the final loss in \Cref{tab:manifold-final-loss}. On the classes of triangulations used for training (spheres, random manifolds, random gluings, and ``near misses'') the final test loss is 0.0576. Interestingly, we find that the test losses for the specific ``extra'' manifold types are all very low despite being topologies that were not seen by the model during training. This suggests that the model could be learning to characterise some property that is common to manifolds, rather than something specific only to the training data. However, the out of sample test data are all positive cases (manifolds), and are ``similar'' to $S^3$, hence may be more easily identifiable. Consistent with this theory, we find that the test loss for the ``random'' manifolds is higher than average, as is the test loss for ``near misses''. This suggests that the model has not specifically learned to examine the vertex link, and that it can be ``confused'' by complex manifolds or triangulations that look like manifolds almost everywhere.

The test loss of $0.0576$ is a substantial improvement from the baseline of guessing completely at random if an input is a closed triangulation or not ($\ln(2)\approx 0.69$). This is especially noteworthy given the size of the model and the limited training data and compute spent on it. While the model does not perfectly learn to classify manifolds, it seems to largely learn what a manifold ``looks like''. The loss of $0.0088$ on the random gluings indicates that transformers are able to reliably identify triangulations that are very far from being closed.

\begin{table}[htpb]
    \centering
    \begin{tabular}{c|c}
        Triangulation Type & Final Loss \\
        \hline
        $S^3$ & 0.0129 \\
        Random manifold & 0.0755 \\
        Near miss & 0.1331 \\
        Random gluing & 0.0088 \\
        $S^2\times S^1$ & 0.0033 \\
        $\mathbb{R}P^3$ & 0.0753 \\
        $L(7, 1)$ & 0.0249 \\
        $L(7, 2)$ & 0.0159 \\
    \end{tabular}
    \caption{Final test loss for manifold classification on different test sets. The first 4 sets are classes that are explicitly included in training.} 
    \label{tab:manifold-final-loss}
\end{table}

\subsection{Scalar Regression}
\label{subsec:scalar-regression}
We examine a series of objective functions of differing complexity. Performance on these objective functions is an indication of the features that the machine learning architecture can learn. All tasks are done with a single model with $d_{\text{model}}=512$, $n_{\text{head}}=8$, $n_{\text{block}}=8$, and $d_{\text{out}} = 64$, and hence $\approx 25$ million parameters. We ran 16 epochs of training on the 4 million scored samples of 10-tetrahedra triangulations of the single vertex 3-sphere, setting aside $4,000$ samples for testing. 

All models are scored according to mean square error (MSE), and the data is standardised before training. This means that a naive model $f(\operatorname{sig})=0$ has an MSE $=1$. We use a coarse peak learning rate grid of $1\times10^{-5}, 3\times 10^{-5}, 1\times 10^{-4}, 3\times 10^{-4}, $ and $1\times 10^{-3}$ to determine the optimal peak learning rate in each case. Similar to the classification task from \Cref{subsec:manifold-property}, this is a small model with a small dataset and a large number of epochs. Hence, performance would likely be higher on larger models with more data.

\paragraph{Dual Graph Self Loops}
The number of self loops in the dual graph is a simple objective function to learn. It contains limited information on the triangulation, and only requires learning to decode the gluing table from the dehydration signature. We find the optimal learning rate to be $1\times10^{-4}$, and observe the loss to converge to zero. The training run converged after two epochs, when we terminated the experiment.

\paragraph{Edge Degree Count}
We examine the edge degree count of the triangulations, that is, the number of edges that have a given degree. 

The set of edge degrees of a triangulation is an interesting and intricate combinatorial property. It allows conclusions on how much a triangulation can be simplified \cite{Jaco09LensSpaces}, it is related to curvature computation in Regge calculus \cite{Regge1961}, and it is used in the assignment of a hyperbolic geometry to a triangulation (Thurston's Gluing Equations) \cite{NeumannZagier1985}. 

We perform our experiments for edge degree = $1, 2, 3, 4, 5$. We choose to stop here because, as a consequence of Euler characteristic computations, the average edge degree of a closed triangulation with $N_{\text{tet}}$ tetrahedra and $v$ vertices is $\frac{6 N_{\text{tet}}}{N_{\text{tet}}+v} < 6$. Hence, knowing about the edges of low degrees $<6$ already determines much about the edges of higher degree. Moreover, as the degree increases, it becomes increasingly challenging to predict the number of edges of such degrees. This is because, whenever an edge has degree $k$, the dual graph must have a closed walk of length $k$ (the converse is not true). Hence, the edge degree count is related to cycle detection, which is a known challenge for deep learning models, both for graph neural networks \cite{chen2020can,loukas2020what,zhang2023rethinking} and for transformers operating on graph inputs \cite{sanford2024understanding,sanford2024transformers}.

We detail the outcome for each edge degree count in \Cref{tab:edge-degree-count}, including the best learning rate selected from the grid and the minimum test loss. We observe that the model learns to essentially perfectly determine the degree-one and degree-two edge counts. For these tasks, we stop training after two epochs. This is mostly unsurprising as these are associated (but not equivalent) to loops and double edges in the dual graph. The loss increases for learning the counts of degree $3$, $4$, and $5$ edges. Nonetheless, in all cases the loss does not exceed $0.51$ indicating a $49\%$ reduction in MSE.

\begin{table}[htpb]
    \centering
    \begin{tabular}{c|c|c}
        Degree & Optimal LR & Minimum loss (test) \\
        \hline
        1 & $3\times10^{-4}$ & $2.0\times 10^{-4}$ \\
        2 & $3\times10^{-4}$ & $0.031$ \\
        3 & $3\times10^{-4}$ & $0.16$ \\
        4 & $1\times10^{-4}$ & $0.23$ \\
        5 & $1\times10^{-4}$ & $0.51$ \\
    \end{tabular}
    \caption{Performance metrics for scalar regression of the count of edges with a given degree.}
    \label{tab:edge-degree-count}
\end{table}

\paragraph{Edge Degree Variance}
The edge degree variance, that is, the mean squared deviation between edge degrees and the average edge degree of a triangulation, is studied for a few reasons. On the one hand, it is a useful aggregate measure for the count of edge degrees above. On the other hand, it is a measure that is known to be highly correlated with some notions of ``complicatedness'' of triangulations, and related to how ``wild'' the curvature of the geometric realisation of the triangulation is \cite{Paixao}.

We find the best learning rate to be $1\times 10^{-4}$ which achieves a minimum test loss of $0.36$. While the model does not learn to perfectly determine the edge degree variance, this is a $64\%$ reduction in MSE, suggesting that the model has at least learned to identify something that is correlated with edge degree variance.

\subsection{Autoregression}
\label{subsec:auto-regression}
We now consider the main goal of the paper: generating triangulations. We conduct three experiments. First, we examine small models on triangulations of different sizes to establish how performance differs on different numbers of tetrahedra. We then examine the scaling law behaviour of the learning task. Finally, we examine the characteristics of the generated data.

\paragraph{Generation Efficiency for Different Triangulation Sizes}
We begin our analysis by performing a small investigation into the generation efficiency for different numbers of tetrahedra. We choose a consistent model architecture of $d_{\text{model}}=512$, $n_{\text{head}}=8$, and $n_{\text{block}}=8$. This model has $\approx 25$ million parameters, and we run 64 epochs of training on 1 million samples for $10 \leq N_{\text{tet}}\leq 20$, leading to approximately 2.5 samples per parameter, consistent across each value of $N_{\text{tet}}$. We leave $4,000$ samples out for holdout. We run a learning rate grid of $10^{-4}, 3\times10^{-4}, 10^{-3}$, and $3\times 10^{-3}$. After training, we generate $10,000$ samples from the best model. For each generated signature we check if it is an acceptable string, if it is valid, if the triangulation it represents is a closed manifold, and if it is a sphere. These results are included in \Cref{tab:ntet-test-results}. We find in all cases, every string generated is unique and represents a combinatorially unique triangulation. In all cases, every generated string is distinct from the training data.

We observe that even for the small models we train, the performance metrics are much better than those established in \Cref{sec:baseline}. Additionally, the percentage of signatures that load, and are valid is particularly high. We additionally observe an approximate exponential decline in the closed and sphere generation efficiency as $N_{\text{tet}}$ increases.

\begin{table}[htpb]
    \centering
    \begin{tabular}{c|c|c|c|c|c|c}
        $N_{\text{tet}}$ & Signature Length & Optimal LR & Acceptable (\%) & Valid (\%) & Closed (\%) & Sphere (\%) \\
        \hline
        10 & 29 & $1\times 10^{-3}$ & 99.55 & 99.53 & 9.95 & 8.08 \\
        11 & 31 & $1\times 10^{-3}$ & 98.32 & 98.32 & 6.2 & 4.89 \\
        12 & 33 & $1\times 10^{-3}$ & 96.02 & 96.01 & 4.03 & 3.08 \\
        13 & 37 & $1\times 10^{-3}$ & 96.64 & 96.64 & 2.35 & 1.83 \\
        14 & 39 & $1\times 10^{-3}$ & 95.87 & 95.86 & 1.80 & 1.25 \\
        15 & 41 & $1\times 10^{-3}$ & 91.29 & 91.28 & 1.22 & 0.92 \\
        16 & 43 & $1\times 10^{-3}$ & 88.53 & 88.52 & 0.71 & 0.47 \\
        17 & 47 & $1\times 10^{-3}$ & 86.02 & 86.00 & 0.42 & 0.33 \\
        18 & 49 & $1\times 10^{-3}$ & 85.59 & 85.56 & 0.37 & 0.20 \\
        19 & 51 & $1\times 10^{-3}$ & 74.73 & 74.72 & 0.12 & 0.07 \\
        20 & 53 & $1\times 10^{-3}$ & 67.21 & 67.16 & 0.05 & 0.05
    \end{tabular}
    \caption{Generation efficiency of the 25 million parameter model for triangulations with different numbers of tetrahedra. Based on 10,000 samples. The signature length does not include the appended \texttt{[SOS]} or \texttt{[EOS]} tokens.}
    \label{tab:ntet-test-results}
\end{table}

\paragraph{IsoFLOP Neural Scaling Laws}
In this section we investigate how the training loss of our models relates to compute, model size, and training data. For this purpose, we estimate the {\em scaling laws} of our problem, with our attention focused on triangulations with $15$ tetrahedra. We focus on this size because $15$ tetrahedra is small enough to expect reasonable performance (based on the above results), and fast computation; but large enough that it is outside of the scope of existing censuses, containing unexplored territory. 

In order to determine the scaling behaviour, we use the IsoFLOP method (i.e. fixed compute method) described in \cite{Chinchilla}. That is, we fix a compute budget $C$, and fit models of a range of different parameter counts (and data sizes) to estimate the optimal achievable loss for the given compute budget. Repeating this for compute budgets from $6\cdot 10^{15}$ to $3\cdot 10^{18}$ FLOPs allows us to extract the relationship between the model size, data, loss and compute.

To ensure validity and consistency, for each FLOP target, we compute in advance the required number of training iterations via the formula $\text{FLOPs}\approx 6ND$, where $N$ is the number of parameters and $D$ is the number of tokens in the training data. While this is just an approximation for the compute budget, it is consistent and independent of hardware and the compiler. Computing the number of iterations in advance is critical for setting up the learning rate schedule, as the model is not fully fit until the schedule is complete. For each FLOP target, the model is completely retrained with its target iterations to ensure the learning rate schedule is complete by the time the final loss is extracted. We perform 63 training runs, detailed in \Cref{appendix:IsoFLOP-runs}. Each model has a consistent aspect ratio ($d_{\text{model}} / n_{\text{block}}$) of $64$ which has been chosen based on investigations outlined in \Cref{appendix:aspect-ratio-tuning}. Because training a large model is computationally expensive, we run a learning rate grid search on a narrow model and scale the learning rate with $\mu P$ \cite{yang2022tensorprogramsvtuning}. We find that for the range of depths being considered, using a consistent \textit{base} learning rate across all the model depths is stable, and results in an acceptable performance loss (performance loss on the order of the noise). The details of this learning rate search are included in \Cref{appendix:peak-learning-rate}.

\begin{figure}[htbp]
    \centering
    \includegraphics[width=0.7\linewidth]{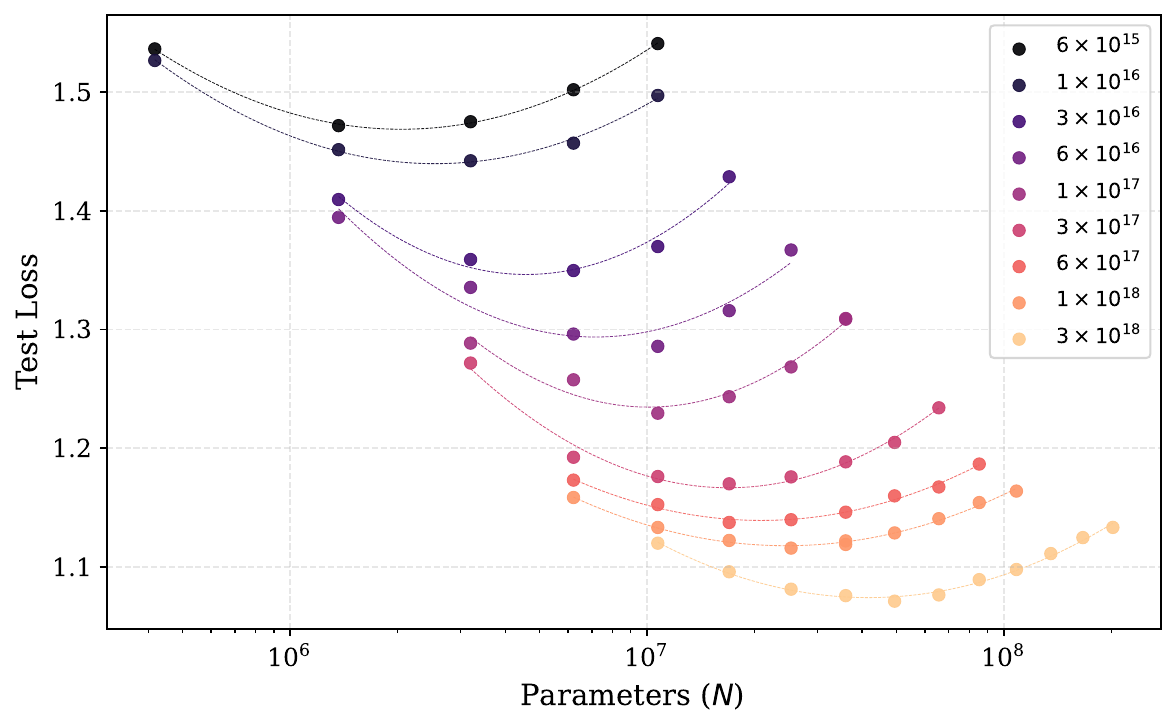}
    \caption{IsoFLOP curves fitted with quadratics. Log $x$-axis}
    \label{fig:IsoFLOP-curves}
\end{figure}

For each fixed compute budget we fit a quadratic to the loss as a function of the log number of parameters, shown in \Cref{fig:IsoFLOP-curves}. For each quadratic fit, we are interested in determining the confidence intervals. However, the residual standard error, $\sigma(C)$, is heteroscedastic (has non-constant variance), ranging from $0.0013$ to $0.015$ with a peak at an intermediate compute. This heteroscedasticity is not captured by a linear trend in $\ln(C)$ ($r=-0.2$). Furthermore, the residuals are heavy-tailed (Shapiro–Wilk $p=0.0004$). Therefore, we use the {\em wild bootstrap} technique on leverage-corrected residuals, which is valid under heteroscedasticity, as our primary method of determining the confidence intervals. As a robustness check, we use a {\em residual bootstrap} from the pooled leverage-corrected residuals, which assumes homoscedastic noise and is conservative here.

\begin{remark}
    The residual standard errors above are consistent with \Cref{appendix:peak-learning-rate}, which had $\sigma=0.004$.
\end{remark}

\begin{remark}
    Although a shared random key was used during initialisation, a permutation test found no correlation between the residuals of the same architecture at neighbouring compute budgets ($\rho=0.02$, $p=0.90$).
\end{remark}

\begin{figure}[htbp]
    \centering
    \begin{subfigure}[b]{0.48\linewidth}
        \centering
        \includegraphics[width=\linewidth]{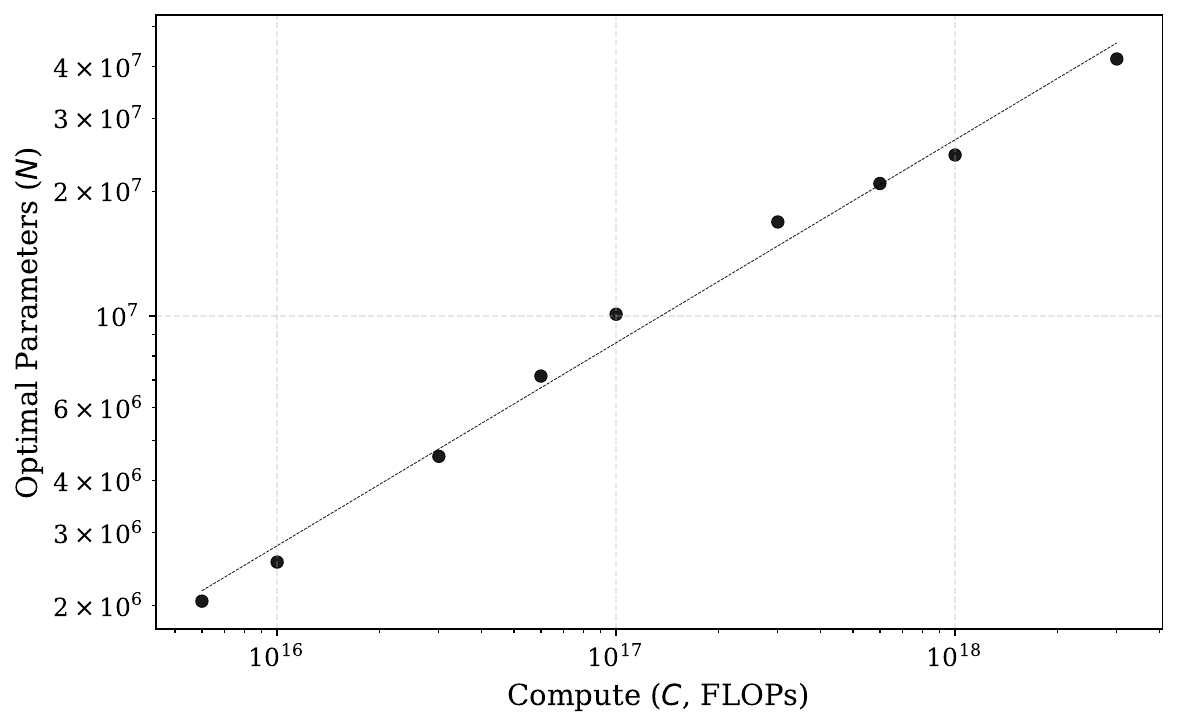}
        \caption{Optimal parameters extracted from quadratic fits vs compute budget. Log-log plot with fitted power law curve.}
        \label{fig:fitted-opt-params}
    \end{subfigure}
    \hfill
    \begin{subfigure}[b]{0.48\linewidth}
        \centering
        \includegraphics[width=\linewidth]{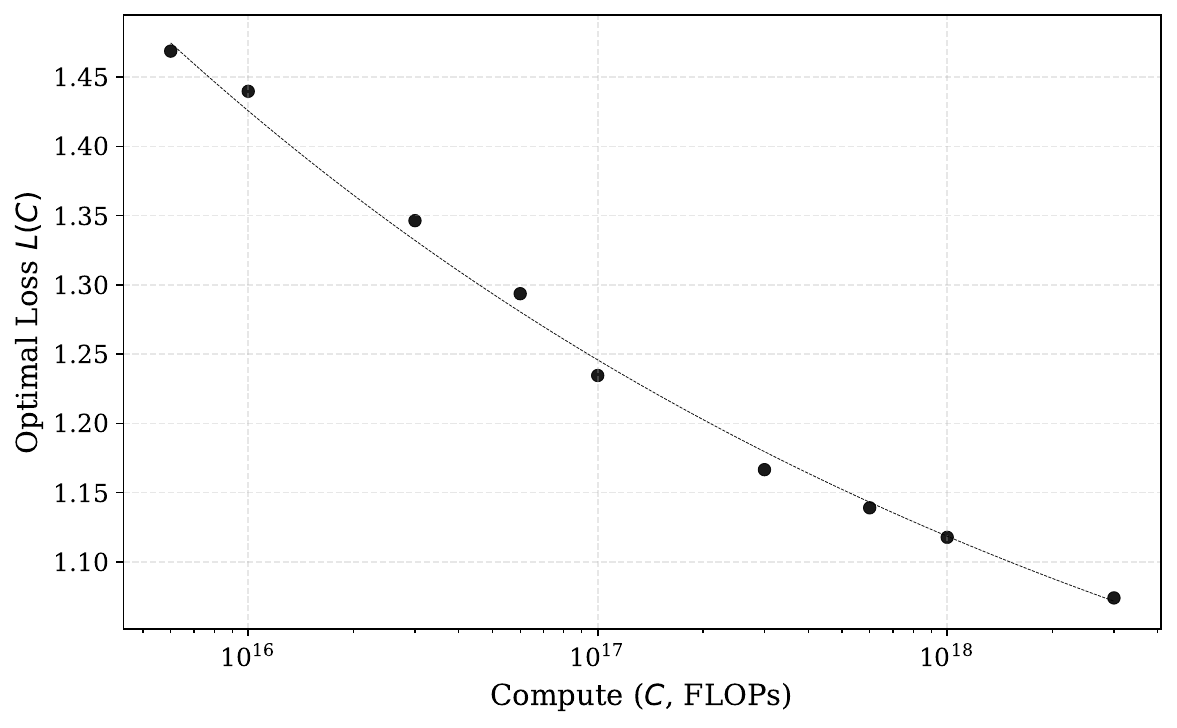}
        \caption{Optimal loss $L(C)$ extracted from IsoFLOP curves vs compute budget $C$. Log $x$-axis plot with fitted power law with constant offset.}
        \label{fig:fitted-opt-loss}
    \end{subfigure}
    \caption{Analysis of parameter and loss scaling.}
    \label{fig:opt-loss-and-params}
\end{figure}

With bootstrapping ($B_{\text{boot}}=1{,}000$), we refit the quadratic to each IsoFLOP test, and extract the optimal parameter count, and the associated loss. We fit a power law $N \sim C^\alpha$ to the optimal parameter counts and find $\alpha = 0.49$, $95\%$ CI $[0.48, 0.50]$ (wild bootstrap), $[0.47, 0.51]$ (residual bootstrap); this fit is shown in \Cref{fig:fitted-opt-params}. The optimal amount of data follows the same parametric form with exponent $1-\alpha$, indicating that data and parameters scale at approximately equal rates. At a reference compute of $1{,}000$ petaFLOPs, the ratio of optimal data to optimal parameters is $234$ tokens per parameter, or $5.6$ samples (that is, dehydration signatures) per parameter. 

\Cref{fig:fitted-opt-loss} shows optimal loss $L(C)$ as a function of compute $C$. We fit the standard form $L(C) = E + aC^{-b}$. Because the residuals contain no isolated outliers, we use least squares loss, and run a grid search on $b$ before refining via Levenberg-Marquardt \cite{levenberg1944method, marquardt1963algorithm}. From this, we find that $E = 0.79$, $95\%$ CI $[0.71,0.86]$ (wild bootstrap) $[0.70, 0.85]$ (residual bootstrap) and $b=0.15$, $95\%$ CI $[0.13, 0.17]$ (wild bootstrap) $[0.12, 0.17]$ (residual bootstrap).

We corroborate these results by separately fitting $L(N, D) = E + AN^{-\alpha'} + BD^{-\beta'}$ following Method 3 of \cite{Chinchilla}. From this, we derive point estimates: a parameter-scaling exponent of $\alpha = \frac{\beta'}{\alpha'+\beta'}=0.50, [0.44, 0.55]$ (wild bootstrap), a data-scaling exponent of $1-\alpha$, a reducible-loss exponent of $b=\frac{\alpha'\beta'}{\alpha'+\beta'}=0.15, [0.13, 0.17]$ (wild bootstrap), and an irreducible loss of $E = 0.81, [0.74, 0.87]$ (wild bootstrap). These estimates are consistent with the previous analysis.

To assess systematic error, we repeat the analysis with cubic rather than quadratic fits to each IsoFLOP curve, and with the lowest and highest compute budgets excluded in turn. Together with the Method 3 estimates, this gives five point estimates for each quantity. We take half their range as the systematic uncertainty, giving $\alpha = 0.49, [0.48, 0.50]$ (stat) $\pm 0.02$ (syst), $E = 0.79, [0.71, 0.86]$ (stat) $\pm 0.10$ (syst), and $b = 0.15, [0.13, 0.17]$ (stat) $\pm 0.04$ (syst). Our estimate for $\alpha$ is robust under both systematic and statistical noise, being of comparable scales. $E$ and $b$ are only weakly constrained in this compute range.


\paragraph{Exploring the Characteristics of the Generated Samples.}
For each compute budget, we select the best model and generate $1{,}000{,}000$ samples. In each case, every acceptable sample represents a combinatorially unique triangulation. In all cases, every generated string is distinct from the training data, indicating that the models are not simply memorising.

\begin{figure}[htbp]
    \centering
    \begin{subfigure}[b]{0.48\linewidth}
        \centering
        \includegraphics[width=\linewidth]{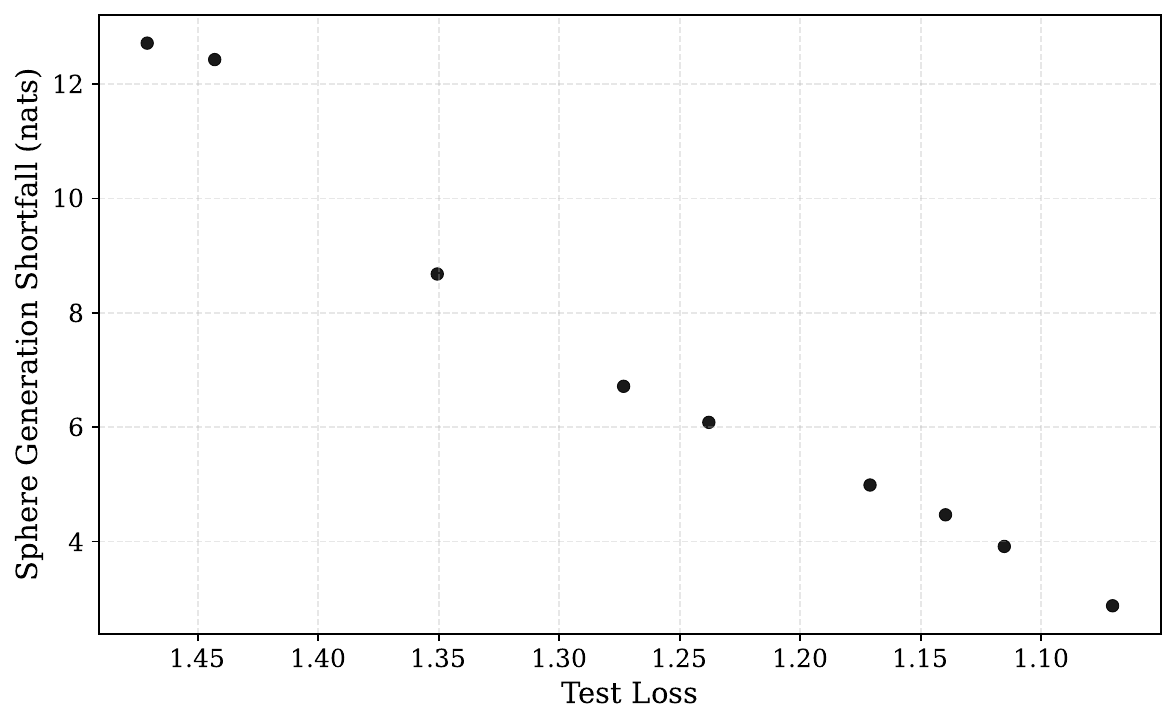}
        \caption{Shortfall (nats) of samples generated that are spheres ($-\ln(p_{\text{sphere}})$) as a function of loss. Inverted $x$-axis.}
        \label{fig:sphere-gen-frac}
    \end{subfigure}
    \hfill
    \begin{subfigure}[b]{0.48\linewidth}
        \centering
        \includegraphics[width=\linewidth]{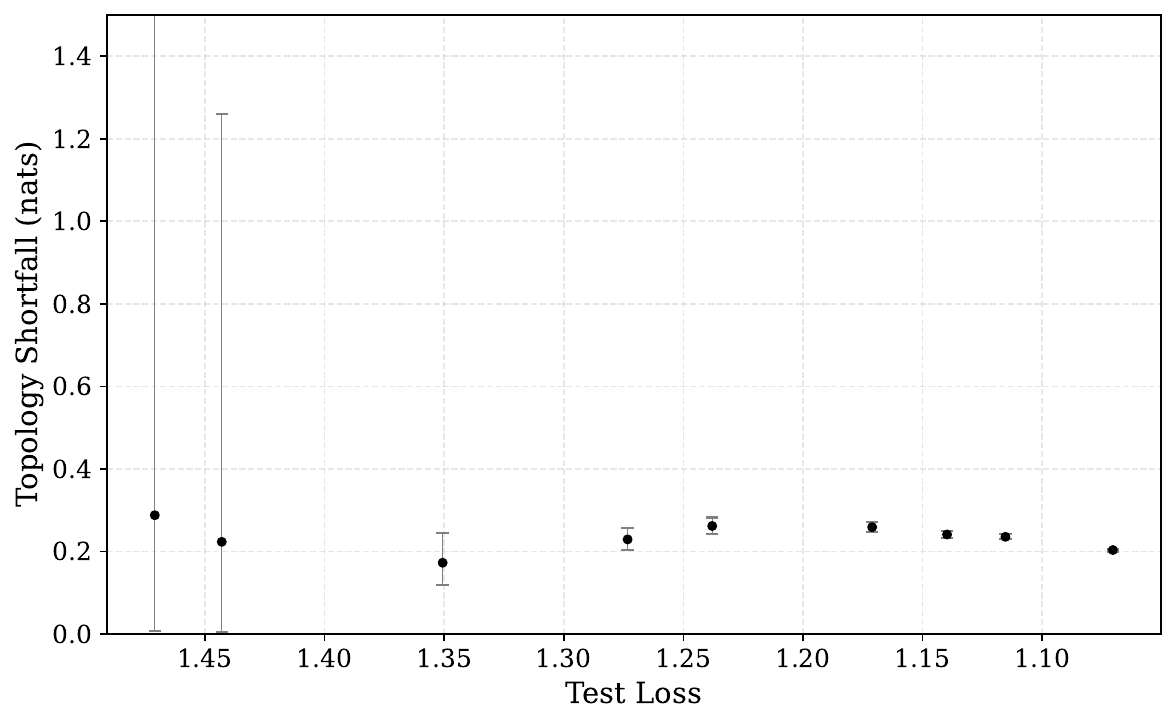}
        \caption{Incremental shortfall of spheres from closed manifolds. Inverted $x$-axis }
        \label{fig:new-closed-frac}
    \end{subfigure}
    \begin{subfigure}[b]{0.48\linewidth}
        \centering
        \includegraphics[width=\linewidth]{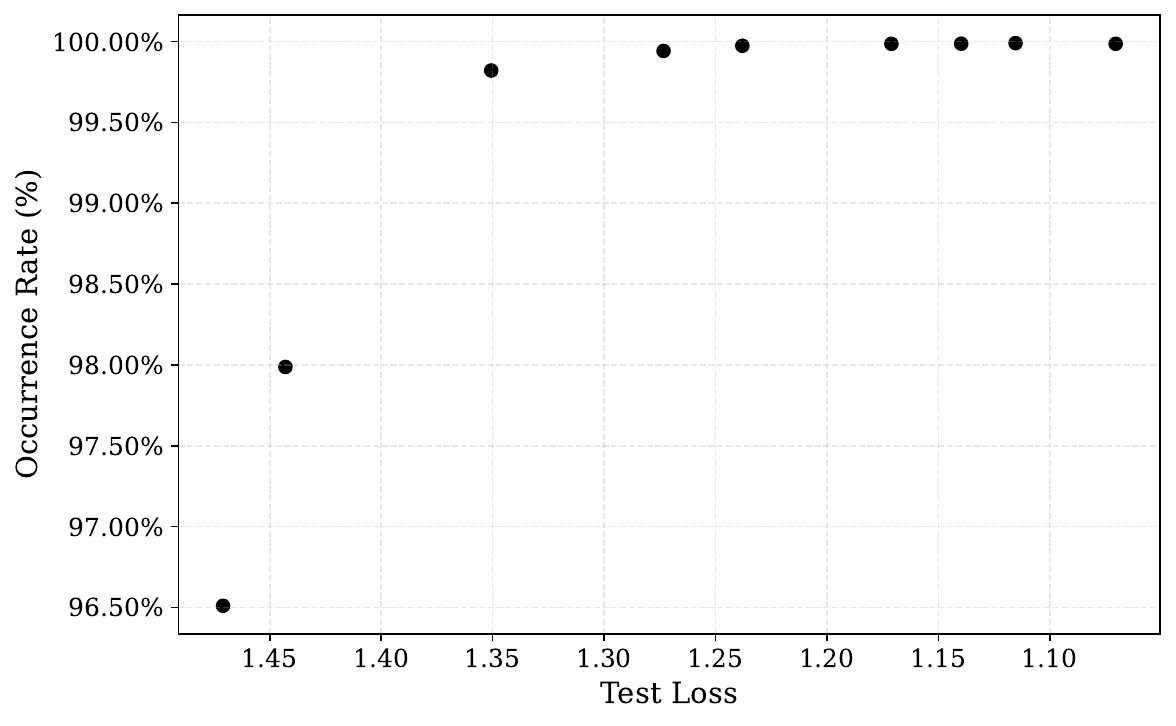}
        \caption{Percentage of vertex links that are orientable vs final test loss.}
        \label{fig:orientable}
    \end{subfigure}
    \hfill
    \begin{subfigure}[b]{0.48\linewidth}
        \centering
        \includegraphics[width=\linewidth]{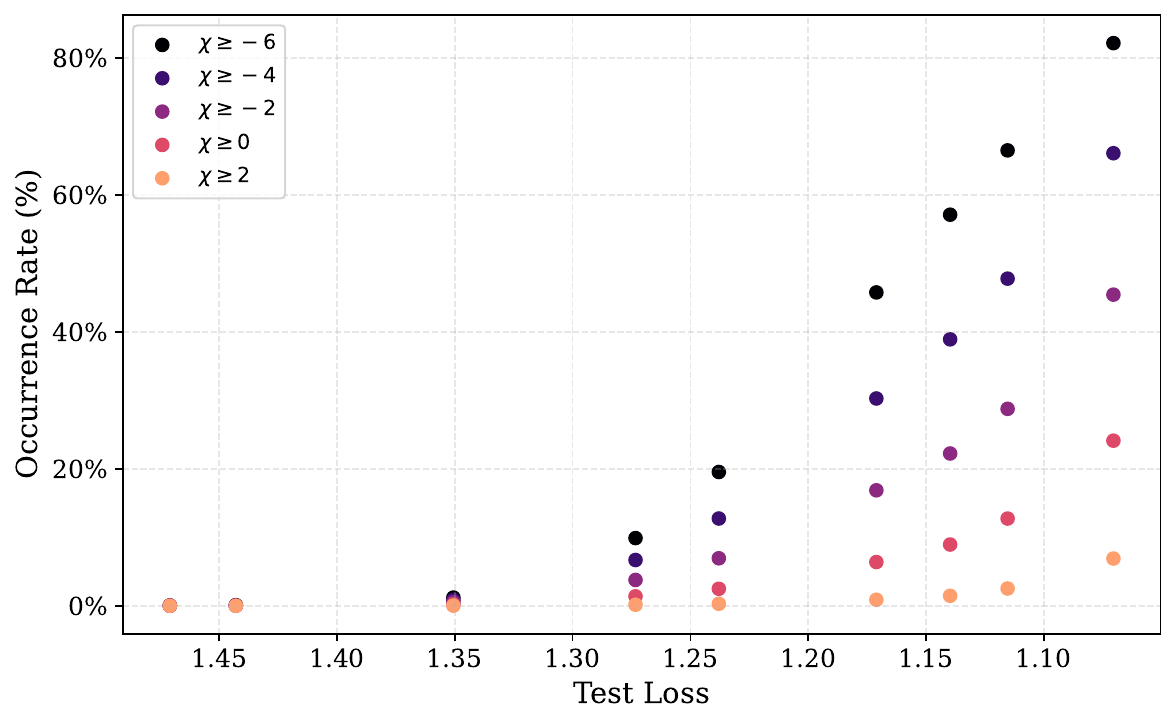}
        \caption{Percentage of generated samples with vertex averaged Euler characteristic $\chi$ exceeding a given threshold.}
        \label{fig:ec-threshold}
    \end{subfigure}
    \caption{Analysis of generated manifolds.}
    \label{fig:gen-manifold-analysis}
\end{figure}

From these samples we compute the sphere generation shortfall in {\em nats}, $-\ln (p_{\text{sphere}})$, where $p_{\text{sphere}}$ is the fraction of generated samples that are spheres. This is shown in \Cref{fig:sphere-gen-frac}. The highest observed sphere rate is $5.63\%$ $[5.59\%, 5.68\%]$ or $2.88$ $[2.87, 2.89]$ nats, achieved by $d_{\text{model}} = 640$. This is significantly higher than the random-gluing baseline of \Cref{sec:baseline}. Pooling the two highest runs and comparing to the best run, the shortfall declines from $12.56$ to $2.88$ nats across the sweep, corresponding to an approximately $16{,}000$-fold increase in the sphere generation rate.

To examine how much of this shortfall improvement is from topology compared to closedness, \Cref{fig:new-closed-frac} shows the conditional shortfall, $-\ln (P(\text{sphere} \;|\; \text{closed}))$. Across the models with more than $100$ closed samples, the conditional shortfall varies between $0.17$ and $0.26$ nats ($77\%$--$84\%$ of closed samples are spheres), with no monotone trend. This is about $1\%$ of the overall $9.68$-nat spread in sphere generation shortfall. Our best model ($d_{\text{model}} = 640$) generates non-sphere closed manifolds at a rate of $1.27\%$ $[1.25\%, 1.29\%]$.

For each valid triangulation, we examine the vertex link. We look at both orientability, and the Euler characteristic. As a baseline we perform random gluings as in \Cref{sec:baseline} on $N_{\text{tet}}=15$ tetrahedra. Including invalid triangulations, we find that 99.6\% of samples generated have non-orientable vertex links, and that the Euler characteristic (averaged across vertices) is $-25$ [$95\%$ central range: $-29$, $-18$]. \Cref{fig:orientable} shows the percentage of vertex links that are orientable against test loss, demonstrating that the models learn the combinatorial property of orientability for the vertex link very quickly. \Cref{fig:ec-threshold} shows the number of samples that have vertex link Euler characteristic greater than some threshold vs the test loss, showing how, as loss improves, the triangulations become steadily ``closer'' to becoming closed manifolds, but that they still typically have negative Euler characteristic. In either case, the transformer generates triangulations that are closer to satisfying the closed-manifold condition.

\begin{remark}
    The vertex link of the triangulations may be an invalid surface; following Regina \cite{regina}, we still define its Euler characteristic as $v-e+f$. Hence, a vertex link Euler characteristic of $-29$ is possible with just $N_{\text{tet}}=15$.
\end{remark}

We pool the above samples together, removing duplicates. This yields $24{,}148$ non-sphere closed manifold triangulations. For each of the non-sphere closed manifold triangulations, using {\em Regina} \cite{regina}, we simplify the triangulation to find $335$ unique ``base'' triangulations. This may contain some duplicate manifolds where a triangulation was unable to be simplified. The prime summands of these triangulations consist of $S^1\times S^2$, genuine lens spaces ($L(p, q)$ with $p\geq2$), and prism manifolds. That is, with the exception of $S^2\times S^1$, all prime summands are spherical. \Cref{tab:gen-type} shows the generation abundance of each prime summand. Of the unique manifolds generated, 84 are prime, 246 consist of two prime components, and 5 consist of three prime components.

In a separate sample, the largest lens space generated was $L(60, 29)$ which conjecturally has $15$ tetrahedra in its minimal triangulation \cite{Matveev1990}.

\begin{remark}
    The number of unique pooled samples is lower than expected. This is because each run generated its samples with the same random key, which leads to high cross run duplication despite the within run duplication being low. This is not a reflection of anything mathematical.
\end{remark}

\begin{table}[ht]
\centering
\begin{tabular}{c|c|c}
    Prime component & $\pi_1$ & Occurrences \\
    \hline
    $L(p,q)$, $p \geq 2$ & $\mathbb{Z}_p$ & 15{,}717 \\
    $S^1 \times S^2$ & $\mathbb{Z}$ & 9{,}080 \\
    $S^3/Q_8$ & $Q_8$ & 15 \\
    Prism manifold & $Q_{12}$ & 5 \\
    Prism manifold & $\mathbb{Z}_3 \rtimes \mathbb{Z}_8$ & 3 \\
    Prism manifold & $Q_8 \times \mathbb{Z}_3$ & 2 \\
    Prism manifold & $\mathbb{Z}_5 \rtimes \mathbb{Z}_8$ & 2 \\
    Prism manifold & $Q_{16} \times \mathbb{Z}_3$ & 1 \\
    Prism manifold & $\mathbb{Z}_3 \rtimes \mathbb{Z}_{16}$ & 1 \\
    \hline
    \textbf{Total} & & \textbf{24{,}826} \\
\end{tabular}
\caption{Occurrences of each prime summand, counted with multiplicity.}
\label{tab:gen-type}
\end{table}

\subsection{Statistical Properties}
We are interested in assessing the statistical properties of our models, and their feasibility in producing reliable statistical estimates over the space of triangulations. In order to do this, we train a large model on $N_{\text{tet}} = 15$ with $9{,}000$ petaFLOPs based on the scaling laws above, achieving a loss of $1.036$, consistent with the fitted curves. 

\paragraph{Spheres}
Using the algorithm from \cite{Altmann26MCMC}, we generate $100{,}000$ triangulations of the $3$-sphere with $N_{\text{tet}}=15$ tetrahedra. We assume that this is a nearly uniform sample. For varying compute budgets, including $9{,}000$ petaFLOPs, we take the best model and, for each triangulation in the sample, compute its negative log likelihood (NLL), that is, $x \in \mathbb{R}$, if the probability for generating the triangulation is $e^{-x}$. 

\begin{figure}
    \centering
    \includegraphics[width=0.5\linewidth]{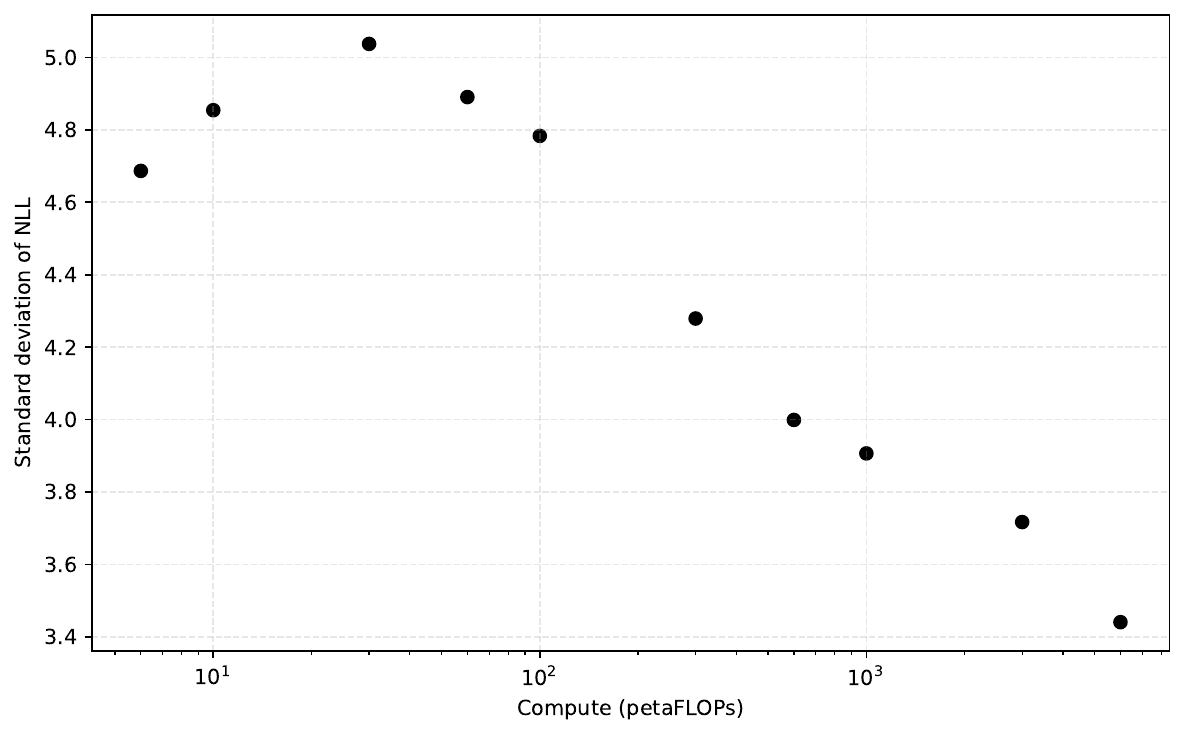}
    \caption{The standard deviation of the NLL vs compute. Log $x$-axis.}
    \label{fig:nll-v-compute}
\end{figure}

We compute the mean, $\mu$, and standard deviation, $\sigma$, of the NLLs. Naturally, lower values of $\sigma$ indicate that the sample is closer to uniform. \Cref{fig:nll-v-compute} shows $\sigma$ as a function of compute. We observe that for larger and larger models, $\sigma$ declines, indicating that the samples are getting increasingly uniform. The smallest observed value of $\sigma$ is $3.44$. However, we also observe a heavy tail in the distribution of the NLLs, with the least likely $3$-sphere triangulation having an NLL of $76.87$. These outliers prevent us from using techniques such as importance sampling to effectively deduce properties of the space of $3$-sphere triangulations from the sample (as done successfully in the MCMC setting in \cite{Altmann26MCMC}).

\paragraph{Other Manifold Types}
We generate and score $12{,}500$ samples of a range of different manifold types. We score each of them with the $9{,}000$ petaFLOP model and compute its mean NLL, $\mu$. This is shown in \Cref{tab:other-manifold-nll}. If the model was generating signatures uniformly across the space of all closed dehydration signatures, regardless of topology, the NLL would be constant. We observe that the mean NLL for $S^3$ is significantly (but not prohibitively) smaller than that of the other manifold types, suggesting that our models are heavily biased towards generating spheres. Additionally, we observe that the mean for lens spaces and $S^2\times S^1$ is significantly lower than for the hyperbolic Weeks manifold or $S^3/Q_8$. 

\begin{table}[htb]
    \centering
    \begin{tabular}{c|c}
        Manifold & $\mu$ \\
        \hline
        $S^3$ & $44.54$ \\
        $\mathbb{R}P^3$ & $50.74$ \\
        $L(5,2)$ & $50.41$ \\
        $L(7,2)$ & $50.29$ \\
        $S^2\times S^1$ & $51.46$ \\
        $S^3/Q_8$ & $54.02$ \\
        \textit{Weeks} & $54.32$
    \end{tabular}
    \caption{Mean negative log likelihood for a range of different manifold types.}
    \label{tab:other-manifold-nll}
\end{table}

\section{Future work and discussion}

Interpretability of our models is currently a hard open question. Instead, we interpret our experimental data. Our work in \Cref{subsec:scalar-regression} and work done in \cite{Costantino_2025} indicate that neural networks are efficient at learning simple combinatorial properties, but struggle to fully learn the details required to distinguish topological type. Additionally, our work in \Cref{subsec:auto-regression} demonstrates that, at least over the compute range we tested, the majority of the improvement in the sphere generation shortfall is due to increased efficiency in generating closed manifolds, not in improving the relative fraction of spheres. While it is not a strict requirement that the model learns to detect topology in order to generate correct samples, the alignment of these two results suggests that there may be a connection. One interpretation could be that in order to learn to produce spherical triangulations, the model needs to be able to learn to detect if a triangulation is a sphere. Despite this, it is important to recognise that the relative rate of spheres compared to closed manifolds is higher than the baseline rate found in \Cref{sec:baseline}, and that the rate is high even at high test losses. Additionally, the non-sphere closed manifolds are all ``simple'' or ``sphere-like'' e.g., there are no hyperbolic manifolds. These observations are not contrary to our hypothesis, but rather suggest that there may be some easily learnable properties that correlate well with being a sphere. 

\paragraph{Increasing diversity of topological types in samples}
The topological types of closed manifolds generated by our $15$-tetrahedron model are listed in \Cref{tab:gen-type}. The list is heavily skewed towards quotients of spheres, but nonetheless exhibits a diversity that is similar to the list of topological types of very small manifolds in the census of closed manifolds shipped with \cite{regina}. We treat this as emerging experimental evidence that, despite the model being trained entirely on sphere triangulations, it is not fully biased to generate spheres. The model can therefore be used as a first stepping stone for exploration over the space of all triangulations. 

There are a number of straightforward extensions of our approach to increase the topological diversity of generated data. First, we can use generated samples from our trained model as seeds to generate new training data to feed into a new model. Iterating this procedure will (a) likely have some effect on the distribution of topological types in samples, and (b) potentially lead to a stable distribution of topological types. This stable distribution, especially if compared to estimates of the true distribution of topological types in the space of all triangulations, is a much stronger indicator for how well transformers pick up on topological properties in the training data.

Another approach is to start with a larger variety of topological types in the training data by, for instance, combining $10^5$ samples for each of the $\approx 13,400$ closed manifolds admitting a triangulation with at most $11$ tetrahedra, available from \cite{BurtonData}, into one data set.

\paragraph{Reinforcement learning} 
Using one of our models as a pretrained baseline, a reasonable extension of our work is to fine-tune this model using reinforcement learning techniques. Given that our pretrained models already have high generation efficiency for closed triangulations (see \Cref{tab:ntet-test-results}), and that vertex link Euler characteristic could be an appropriate candidate for a scoring function, we see this as a promising approach to improve our models.  

\paragraph{Other model architectures}
Autoregression can have difficulty on long sequences as it is unable to correct chance mistakes during generation. Diffusion-style models are a possible approach to address this, and investigation into this could yield improved generation efficiencies.

Additionally, we have only looked into transformers of a very specific type as a model architecture. However, triangulations have a natural graph structure given by their dual graphs. Derivatives of graph neural networks may be good at exploiting this dual-graph structure.

\paragraph{Machine-learning-friendly encoding of triangulations}
In this article, we use dehydration signatures to represent triangulations. While this works well with our GPT-2-style transformers, there are some natural limitations to this framework. Most importantly, dehydration signatures are only defined for triangulations up to 25 tetrahedra. While extending dehydration signatures to larger triangulations is straightforward, other tweaks may be possible to expose the combinatorics of a triangulation even more efficiently to a transformer. 

A systematic optimisation of triangulation encodings for machine learning applications has the potential to substantially affect training success: initial training runs with {\em Regina} isomorphism signatures showed much lower generation efficiency compared to dehydration signatures in similar models.

\bibliographystyle{plain}
\bibliography{bibliography}

\newpage
\begin{appendices}
\section{Estimating Unique Dehydration Signatures}
\label{appendix:estimating-dehydration-count}
To estimate the number of unique, acceptable dehydration signatures, we construct a tree where each proper leaf of the tree represents a unique dehydration signature. Specifically, for a fixed $N_{\text{tet}}$, we break each dehydration signature into 4 parts: A binary string encoding the ``new tetrahedra'' information, a series of letters encoding the ``where'' information, a list of numbers from 0--3 representing which face is glued in each ``explicit'' gluing, and a list of permutations of the vertices (now in $S_3$ as the face data is separated). Because any choice of vertex permutation for a fixed face leads to an acceptable dehydration signature, we do not include it in our tree, multiplying by $6^{N_{\text{tet}}+1}$ at the end of our estimation. Additionally, to form a dehydration signature, we need to have a multiple of 8 bits in the gluing bit string, these remaining bits are essentially ``padding'' and can be set to any value without impacting the decoding process, leading to a further multiplicative factor that we include at the end.

To construct the tree, we start at the root where the binary string is empty, the ``where'' string is empty, and the face gluing list is empty. At any node we can choose to add a new tetrahedron by concatenating a $1$ to the binary string, or we can choose to perform a gluing to an ``explicit'' tetrahedron, concatenating a $0$ to the binary string, and adding the letter of the target tetrahedron and face id to their respective data. When we say ``explicit'', we allow forward references where we glue to a tetrahedron that hasn't been declared in the binary string yet as these \textit{are} acceptable. Before listing these children, we remove any children that are not \textit{partially valid} by checking if they lead to any obvious double gluings or gluings of a face to itself. Because the partial check is not perfect, the tree contains \textit{dead ends} in it i.e., nodes that have no partially valid children but are not themselves complete.

To form our estimate, we use random walks as introduced in \cite{knuth1975estimating}, with non-uniform selection probabilities. That is, we start at the root of the tree, randomly select a child node according to some probability distribution, and repeat down the tree. Each walk gives an unbiased estimate of the number of leaves, $\hat{N} = 1/\prod_i p_i$, where $p_i$ is the probability of the child chosen at step
$i$, and $\hat{N}=0$ if the walk fails to reach an acceptable dehydration signature. Specifically, on the dehydration signature tree, if a child that involves gluing a new tetrahedron exists, it gets a weight of $1$ if it is the only child, and otherwise a weight of exactly $p_{\text{new}}=\frac{n_1}{n_0+n_1}$ where $n_0 = N_{\text{tet}} + 1 - \#\{\text{zeros used}\}$ and $n_1 = N_{\text{tet}} - 1 - \#\{\text{ones used}\}$ where the number of ones and zeros used refers to the binary string. The remaining children get the remaining weight divided evenly between them. With this method, in the largest case of $N_{\text{tet}}=20$ there are still $\approx 5\%$ of samples that are not dead ends.

For each $N_{\text{tet}}$ we perform $10{,}000$ walks as above and compute the average of $\hat{N}$ over all walks. This is our estimate of the number of dehydration signatures. It is important to acknowledge that, despite this being an unbiased estimate, the variance can be extremely high and ``heavy tailed'', and it is difficult to estimate what it is. Hence, to provide some diagnostic confidence, we plot the running average and visually confirm our estimate has converged by checking for an absence / stabilisation of ``jumps''. This is not a perfect diagnostic.

\section{IsoFLOP Training Run Details}
\Cref{tab:IsoFLOP} includes the model shape, FLOPs, and data used for each run of the IsoFLOP scaling law test.

\label{appendix:IsoFLOP-runs}
\begin{table}[htpb]
    \centering
    \footnotesize
    \setlength{\tabcolsep}{3pt}
    \begin{tabular}{rrr *{9}{r}}
        \toprule
        & & & \multicolumn{9}{c}{Target petaFLOPs} \\
        \cmidrule(l){4-12}
        $d_{\text{model}}$ & $n_{\text{block}}$ & $N$ & 6 & 10 & 30 & 60 & 100 & 300 & 600 & 1\,000 & 3\,000 \\
        \midrule
        128  & 2  & 417\,280     & 58 & 97 & -- & --  & --  & --  & --  & --  & --   \\
        192  & 3  & 1\,365\,696   & 18 & 30 & 89 & 179 & --  & --  & --  & --  & --   \\
        256  & 4  & 3\,200\,512   & 8  & 13 & 38 & 76  & 127 & 381 & --  & --  & --   \\
        320  & 5  & 6\,216\,640   & 4  & 7  & 20 & 39  & 65  & 196 & 392 & 654 & --   \\
        384  & 6  & 10\,708\,992  & 2  & 4  & 11 & 23  & 38  & 114 & 228 & 380 & 1139 \\
        448  & 7  & 16\,972\,480  & -- & -- & 7  & 14  & 24  & 72  & 144 & 240 & 719  \\
        512  & 8  & 25\,302\,016  & -- & -- & -- & 10  & 16  & 48  & 96  & 161 & 482  \\
        576  & 9  & 35\,992\,512  & -- & -- & -- & --  & 11  & 34  & 68  & 113 & 339  \\
        640  & 10 & 49\,338\,880  & -- & -- & -- & --  & --  & 25  & 49  & 82  & 247  \\
        704  & 11 & 65\,636\,032  & -- & -- & -- & --  & --  & 19  & 37  & 62  & 186  \\
        768  & 12 & 85\,178\,880  & -- & -- & -- & --  & --  & --  & 29  & 48  & 143  \\
        832  & 13 & 108\,262\,336 & -- & -- & -- & --  & --  & --  & --  & 38  & 113  \\
        896  & 14 & 135\,181\,312 & -- & -- & -- & --  & --  & --  & --  & --  & 90   \\
        960  & 15 & 166\,230\,720 & -- & -- & -- & --  & --  & --  & --  & --  & 73   \\
        1024 & 16 & 201\,705\,472 & -- & -- & -- & --  & --  & --  & --  & --  & 60   \\
        \bottomrule
    \end{tabular}
    \caption{Number of training samples ($\times 10^6$) for each IsoFLOP run, by model size and target compute. $N$ is the exact parameter count. Each model has $n_{\text{head}} = n_{\text{block}} = d_{\text{model}}/64$. ``--'' indicates configurations not run at that budget.}
    \label{tab:IsoFLOP}
\end{table}

\section{Aspect Ratio Tuning}
\label{appendix:aspect-ratio-tuning}
We examine the influence of the aspect ratio of the model on the final test loss. We do not perform an exhaustive study on this topic, running only small tests of $10^{16}$ FLOPs on a target parameter count of 6.2 million, specifically on the $15$-tetrahedron dataset at one epoch on a subset of the $1.2$ billion sphere dataset. This dataset is chosen because it is the testing ground for the IsoFLOP scaling law work. With 41 tokens per signature, this is approximately 1 sample per parameter.

We run a range of different aspect ratios, attempting to find models with as close to 6.2 million parameters as possible. However, there is some variance resulting from the need to keep the depth an integer, and the width a multiple of 32. Hence we use (width, depth)-pairs of $(512, 2)$, $(416, 3)$, $(352, 4)$, $(320, 5)$, $(288, 6)$, $(256, 8)$, and $(224, 10)$. For each model we fix the head size at $32$ for consistency. We train each model for approximately $10^{16}$ FLOPs to account for the minor differences in the parameter counts.

We use a coarse learning rate grid of $6\times 10^{-4}, 10^{-3},3\times 10^{-3},6\times 10^{-3},$ and $10^{-2}$. The exact parameter counts, optimal learning rate, and final loss for each model are shown in \Cref{tab:aspect-ratio-tune}. We observe some performance differences at the extremes of the tested aspect ratios, but for aspect ratios between $48$ and $138.7$ we find a range of $0.015$ on the final test loss, aligning with the theory that model performance is robust to minor aspect ratio variation in the vicinity of the optimum. Hence, for subsequent experiments we assume a convenient aspect ratio of 64, acknowledging that performance could be marginally improved by specifically tuning it for each run.

\begin{table}[htpb]
    \centering
    \begin{tabular}{c|c|c|c|c|c}
        Width & Depth & Aspect Ratio & Parameter Count & Optimal learning rate & Final test loss \\
        \hline
        224 & 10 & 22.4 & 6,086,528 & $3\times 10^{-3}$ & 1.473 \\
        256 & 8 & 32 & 6,359,552 & $3\times 10^{-3}$ & 1.444 \\
        288 & 6 & 48 & 6,041,088 & $3\times 10^{-3}$ & 1.386 \\
        320 & 5 & 64 & 6,216,640 & $3\times 10^{-3}$ & 1.378 \\
        352 & 4 & 88 & 6,022,720 & $1\times10^{-3}$ & 1.382 \\
        416 & 3 & 138.7 & 6,313,632 & $3\times 10^{-3}$ & 1.371 \\
        512 & 2 & 256 & 6,387,712 & $3\times 10^{-3}$ & 1.402
    \end{tabular}
    \caption{Performance for different aspect ratios with a target parameter count of 6.2 million. Each model is trained for $10^{16}$ FLOPs (up to minor variance from batching).}
    \label{tab:aspect-ratio-tune}
\end{table}

\section{Peak Learning Rate Tuning}
\label{appendix:peak-learning-rate}
Running a grid search for the optimal learning rate for large models is computationally inefficient. This section focuses on running a search on $d_{\text{model}}=256$ so we can use $\mu P$ to scale the learning rate to larger models. We do this independently for each depth from 2 to 16 (as these are the depths used) with 4 heads, trained for one epoch with a dataset approximately equal to the number of parameters, that is, $\approx 1$ sample per parameter. In keeping with the work of \cite{Chinchilla}, we want to find a stable learning rate that is not excessively low. As such we run 4 copies of the learning rate grid $(0.6, 1, 2, 3, 4, 6)\times 10^{-3}$ with identical datasets but different initialisations. We make two observations from these tests. First, there is a clear ``stable'' and ``unstable'' learning rate range. The centre-point of the stable learning rate range is recorded in \Cref{tab:learning-rate-tune} as ``Optimal learning rate'', along with the number of samples and parameter count. Second, within the ``stable'' learning rate range there is some relation between average final test loss and learning rate. However, the variance of the average loss across different stable learning rates is on the order of the variance between replications, making this a minor factor. Hence, to ensure stability during training, and to reduce confounding effects, we fix the base learning rate at $0.6\times 10^{-3}$ for all depths, noting that using a per depth learning rate would yield lower individual losses, but over the small range of depths we are testing, such a minor optimisation is unlikely to significantly contribute to the outcomes.

\begin{table}[htpb]
    \centering
    \begin{tabular}{c|c|c|c}
        Depth & Parameter count & Number of samples & Optimal learning rate \\
        \hline
        2 & 1,620,992 & 1,572,864 & $4\times 10^{-3}$ \\
        3 & 2,410,752 & 2,359,296 & $2\times 10^{-3}$ \\
        4 & 3,200,512 & 3,145,728 & $3\times 10^{-3}$ \\
        5 & 3,990,272 & 3,932,160 & $1\times 10^{-3}$ \\
        6 & 4,780,032 & 4,718,592 & $3\times 10^{-3}$ \\
        7 & 5,569,792 & 5,505,024 & $3\times 10^{-3}$ \\
        8 & 6,359,552 & 6,291,456 & $2\times 10^{-3}$ \\
        9 & 7,149,312 & 7,077,888 & $2\times 10^{-3}$ \\
        10 & 7,939,072 & 7,864,320 & $3\times 10^{-3}$ \\
        11 & 8,728,832 & 8,650,752 & $1\times 10^{-3}$ \\
        12 & 9,518,592 & 9,437,184 & $2\times 10^{-3}$ \\
        13 & 10,308,352 & 10,223,616 & $1\times 10^{-3}$ \\
        14 & 11,098,112 & 11,010,048 & $1\times 10^{-3}$ \\
        15 & 11,887,872 & 11,796,480 & $0.6\times 10^{-3}$ \\
        16 & 12,677,632 & 12,582,912 & $0.6\times 10^{-3}$
    \end{tabular}
    \caption{Optimal learning rate for different model depths.}
    \label{tab:learning-rate-tune}
\end{table}

\end{appendices}

\end{document}